\documentclass[english]{article}
\usepackage[T1]{fontenc}
\usepackage[latin1]{inputenc}
\usepackage{amsmath}
\usepackage{amssymb}

\makeatletter

\newcommand{\noun}[1]{\textsc{#1}}

 \newcommand{\lyxaddress}[1]{
   \par {\raggedright #1 
   \vspace{1.4em}
   \noindent\par}
 }

\usepackage{babel}
\makeatother

\title{Self-Correction of Transmission\\
on Regular Trees\protect}

\begin{document}

\author{Alberto Gandolfi\protect\footnote{Research partially supported by italian MIUR PRIN Grant \# 2004015228} and Roberto Guenzani\protect\footnote{Research partially supported by italian MIUR PRIN Grant \# 2004015228}}

\maketitle

\begin{abstract}
We consider noisy binary channels on regular trees and introduce periodic
enhancements consisting of locally self-correcting the signal in blocks
without break of the symmetry of the model. We focus on the realistic
class of within-descent self-correction realized by identifying all
descendants $k$ generations down a vertex with their majority. We
show that this also allows reconstruction strictly beyond the critical
distortion. We further identify the limit at which the critical distortions
of within-descent $k$ self-corrected transmission converge, which
turns out to be the critical point for ferromagnetic Ising model on
that tree. We finally discuss how similar phenomena take place with
the biologically more plausible mechanism of eliminating signals which
are locally not coherent with the majority.
\end{abstract}

\newpage

\section{Introduction}

We consider a binary channel on a regular tree, as in \cite{key-1},
with a distortion rate $\varepsilon>0$ at every transmission and
are interested in the reconstruction of the starting bit $\sigma_{0}$
from the signals $\sigma_{W_{n}}$ at the $n$-th generation of the
tree. We focus on the majority rule, by which $\sigma_{0}$ is reconstructed
as the symbol having majority in $\sigma_{W_{n}}$.

In \cite{key-1} it is shown that for regular trees the majority rule
is asymptotically equivalent to the optimal maximum-likelihood rule,
and that there is a critical distortion $\bar{\varepsilon}_{c}=\frac{\sqrt{r}-1}{2\sqrt{r}}$
such that for $\varepsilon>\bar{\varepsilon}_{c}$ no asymptotic reconstruction
takes place and for $\varepsilon<\bar{\varepsilon}_{c}$ there is
asymptotic reconstruction; see also \cite{key-11} for a review and
\cite{key-12} for a dynamical version of these results.

The aim of this paper is to investigate how a non-symmetry breaking
mechanism of correction performed while transmitting the signal can
improve reconstruction by either majority or maximum likelihood. To
this purpose, we propose a local self-correction method by which the
signal is periodically enhanced in blocks formed within the generations.
The enhancement uses majority rule and consists of taking all signals
in a block and changing them to all agree with their majority value
(with random choice to break tie). The self-correction is based on
the information available at the level of interest, and thus can in
principle be performed while the signal is transmitted. From every
vertex the transmission is then continued as it used to be in the
original mechanism and the symmetry of the model is not broken.

It is easy to see that with non-local enhancement one can reconstruct
beyond the critical distortion: in fact, by forcing all vertices of
each generation to agree with their majority, one can reconstruct
for every $\varepsilon\in[0,\frac{1}{2})$. However, such correction
involves taking majority on larger and larger blocks, which is not
an implementable strategy. 

A slightly less expensive self-correction strategy consists of using
blocks of fixed size $M$ (as soon as the generation is large enough)
and then performing self-correction at every generation. In section
2 we show that for any noise level $\varepsilon<\frac{1}{2}$ it is
possible to achieve reconstruction in this way with sufficiently large
block size $M$. This procedure has the advantage of involving only
a bounded number of within generation information exchange in self-correcting
a block, and thus could in principle be implemented by a real machine.
However, it still involves a very large number of within generation
operations, performed at each generation: if the cost of each such
operation is not zero (as in basically all reasonable situations)
then the total cost might become too high.

We, therefore, restrict our attention, in the sequel, to a self-correction
mechanism which contains costs by performing self-correction less
often, and which has the additional advantage of being performed within
the descent of some signal involved in the previous correction. This
within descent self-correction reduces implementation costs, and allows
signals to be dispersed and loose contact after their involvement
in the enhancement, a feature which could be meaningful in a realistic
setting. The within-descent self-correction at level $k$ is performed
by taking each vertex at some $lk$-th generation, $l\in\mathbb{N}$,
considering its $r^{k}$ descendants $k$ generations down, and then
changing them to agree with their majority (randomly breaking ties).

At first sight, it is not even obvious that such reconstruction improves
upon the non self-corrected transmission, but in section 3 we show
that, except for $k=1$ and $r=2$, the within-descent self-correction
at level $k$ strictly increases the critical distortions, and thus
is an effective enhancement. The proof is based on the comparison
between the self-correction based on the majority transformation with
one correction based on random transformation which leaves the critical
points unchanged.

The rest of the work is devoted to identifying the limit of the critical
distortions of the within-descent self-correction of level $k$ as
$k$ diverges. Although it might seem that such mechanism is almost
useless for large $k$, it turns out that instead it improves the
transmission further.

To identify the large $k$ limit, in section 4 we exploit the correspondence
with the Ising model. In fact, it is easy to see that, for regular
trees, the reconstruction problem is equivalent to the free boundary
conditions phase transition of the ferromagnetic Ising model on the
tree with inverse temperature $\beta$ such that $1-2\varepsilon=\tanh(\beta)$.
Such transition occurs at the critical inverse temperature $\bar{\beta}_{c}$
such that for $\beta>\bar{\beta}_{c}$ the free boundary Ising model
is convex combination of the extremal states (see \cite{key-1} for
a detailed description). On the other hand, the Ising model undergoes
its regular phase transition (with boundary conditions) at a lower
inverse temperature $\beta_{c}<\bar{\beta}_{c}$. In terms of $p=\tanh(\beta)$
and on a regular tree with forward branching rate $r$, we have $p_{c}=\tanh(\beta_{c})=\frac{1}{r}$
(as shown originally in \cite{key-8}) and $\bar{p}_{c}=\tanh(\bar{\beta}_{c})=\frac{1}{\sqrt{r}}$
(as shown in \cite{key-2,key-3,key-4}).

Our self-correction at level $k$ introduces thus new critical values
$1-2\varepsilon_{c}(k)=p_{c}(k)=\tanh(\beta_{c}(k))<\bar{p}_{c}$
and our main result is a bound on $p_{c}(k)$ showing that $\lim_{k\rightarrow\infty}p_{c}(k)=p_{c}$,
the regular Ising model phase transition point. Such estimate is derived
by introducing the FK representation of the Ising model and then comparing
the information carried by the FK tree of the origin against the external
{}``noise'' produced by all other freely fluctuating clusters of
vertices. We think that this comparison, which is based on Gaussian
approximation and large deviation techniques, has an interest in itself
as it gives a very natural way of evaluating the information available
on the tree.

In section 5 we remark that the majority self-correction is not biologically
feasible, and introduce, instead, a minority removal self-correction
which consists of self-correcting a generation by removing the elements
not belonging to the majority. Since this leaves at least $r^{k}/2$
descendants, nothing really changes, and such correction also improves
upon normal reconstruction up to the Ising model critical point. As
we discuss, this, however, seems to indicate a peculiar phenomenon:
it looks like that accepting the risk of creating uniform incorrect
regions ({}``tumors'') increases the resistance of inheritance to
distortion. Whether this is a biologically meaningful statement should
be further investigated with many bits models and realistic parameters.

There remain several open issues. First of all, our bounds on $p_{c}(k)$
in section 4 are not sharp. Also, our analysis has been performed
either for correction each $k=1$ steps using large block size $M$
or for correction every $k$ steps with $M=r^{k}$: we do not deal
with the generic case of correcting blocks of size $M$ each $k$
generations. Solving the two issues above would then allow to treat
the main open problem left by the present work: if one is to reconstruct
the signal at a fixed generation $n$ and if within generation transmission
has some given cost, it would be natural to introduce a correspondence
between within generation transmission costs and gain in reconstruction
probability, and then look for the self-correction algorithm with
optimal $k$ and $M$.

\section{Large Block Reconstruction}

We consider regular trees $T^{(r)}$with forward branching rate $r>0$.
The $n$-th level of the tree is indicated by $T_{n}^{(r)}$ and $T_{\rightarrow n}^{(r)}$
represents the tree up to and including the $n$-th level. Vertices
$v$ of $T^{(r)}$ are then identified by coordinates $v=(n,s)$ where
$n$ is the level and $s=1,...,r^{n}$ numbers the vertices at the
same level. Signals or configurations are variables $\{\sigma_{v}\}_{v\in T^{(r)}}$,
$\sigma_{v}\in\{-1,1\}$, and their distribution is specified by taking
$\varepsilon>0$, $P_{\varepsilon}(\sigma_{0}=1)=1/2$ and for each
vertex $v$ and predecessor $^{\leftarrow}v$, $P_{\varepsilon}(\sigma_{v}=\sigma_{^{\leftarrow}v})=1-\varepsilon$
independently of all other pairs. Reconstruction under majority rule
on $(T^{(r)},P_{\varepsilon})$ takes place if

\begin{eqnarray*}
0<\liminf_{n}\Delta_{n}(P_{\varepsilon}) & =: & \liminf_{n}\left(P_{\varepsilon}(S_{n}>0|\sigma_{0}=1)-P_{\varepsilon}(S_{n}<0|\sigma_{0}=1)\right)\\
 & = & \liminf_{n}E_{\varepsilon}|P_{\varepsilon}(\sigma_{0}=1|S_{n})-P_{\varepsilon}(\sigma_{0}=-1|S_{n})|\end{eqnarray*}
\begin{equation}
\end{equation}
where $S_{n}=\sum_{v\in T_{n}^{(r)}}\sigma_{v}$.

We first consider self-correction performed at each step using large
blocks. We fix an integer $M>0$ and let $\tilde{n}=\max\{ k:r^{k}\leq M\}$.
We then consider the $\tilde{n}$-th generation as block $0$, and
partition each of the following generations into blocks of size $M$
as follows: vertices $v=(n,s)\in T_{n}^{(r)}$ are partitioned into
$\left\lfloor \frac{r^{n}}{M}\right\rfloor $ blocks of vertices with
consecutive coordinates $s$, and possibly one block of $r^{n}-\left\lfloor \frac{r^{n}}{M}\right\rfloor $M
vertices, which is from now on discarded without affecting the argument
which follows. Each block $B$ is then connected to all blocks $B'$
such that there are two vertices $v\in B$ and $v'\in B'$ which are
connected on $T^{(r)}$. One can easily see that considering blocks
as renormalized vertices and connections between them as renormalized
bonds we have a new tree $\bar{T}^{(r)}$ with forward branching $r$
at all vertices $\bar{v}\in\bar{T}_{n}^{(r)}$, $n\geq1$, and branching
rate $r_{0}\leq r$ at the starting vertex $\bar{v}_{0}.$ The branching
rate of $\bar{T}^{(r)}$ is thus again $r$. 

Next, we consider self-corrected variables, which are required to
be constant on blocks:

\begin{equation}
\Sigma_{M}=\{\sigma\in\{-1,1\}^{T^{(r)}}\textrm{ such that }\sigma_{v}\textrm{ is constant on each block}\},\end{equation}
and the self-correction map $\Phi_{M}:\{-1,1\}^{T^{(r)}}\rightarrow\Sigma_{M}$
defined by

\[
(\Phi_{M}\sigma)_{v}=\left\{ \begin{array}{ll}
sign(\sum_{v\in B}\sigma_{v}) & \qquad\begin{array}{l}
\textrm{if }v\in B\subseteq T^{(r)}\setminus T_{\rightarrow(\tilde{n}-1)}^{(r)}\\
\textrm{and }\sum_{v\in B}\sigma_{v}\neq0\end{array}\\
\\\\Z & \qquad\begin{array}{l}
\textrm{if }v\in B\subseteq T^{(r)}\setminus T_{\rightarrow(\tilde{n}-1)}^{(r)}\\
\textrm{and }\sum_{v\in B}\sigma_{v}=0\end{array}\\
\\\\\sigma_{v} & \qquad\:\;\textrm{if }v\in T_{\rightarrow(\tilde{n}-1)}^{(r)},\end{array}\right.\]
\begin{equation}
\end{equation}
where $Z\in\{-1,1\}$ is a symmetric random variable. 

The transmission is then self-corrected by the map $\Phi_{M}$at every
step: $\sigma_{T_{\rightarrow(n-1)}^{(r)}}\in\Sigma_{M}$ generates
$\sigma_{T_{\rightarrow n}^{(r)}}\in\{-1,1\}^{T^{(r)}}$ as usual,
and then we take $\Phi_{M}\left(\sigma_{T_{\rightarrow n}^{(r)}}\right)\in\Sigma_{M}$.
The distribution $P_{\varepsilon,M}$ of the self-corrected configuration
is then recursively defined by $P_{\varepsilon,M}\left(\left.\sigma_{T_{\rightarrow n}^{(r)}}\right|\sigma_{T_{\rightarrow(n-1)}^{(r)}}\right)=P_{\varepsilon}\left(\left.\Phi_{M}^{-1}\sigma_{T_{n}^{(r)}}\right|\sigma_{T_{\rightarrow(n-1)}^{(r)}}\right)$. 

We then take configurations on the renormalized tree $\bar{T}^{(r)}$
to be $\bar{\sigma}_{\bar{v}}$ if $\bar{v}$ represents the block
$B$ and $(\Phi_{M}\sigma)_{v}=\bar{\sigma}_{\bar{v}}$ for all $v\in B$,
and indicate by $\Psi_{M}:\Sigma_{M}\rightarrow\bar{\Sigma}_{M}$,
with $\bar{\Sigma}_{M}=\{\bar{\sigma}_{\bar{v}},\bar{v}\in\bar{T}^{(r)}\}=\{-1,1\}^{\bar{T}^{(r)}}$,
the renormalizing transformation. Renormalized configurations are
described by $\bar{P}_{\varepsilon,M}=\Psi_{M}\circ P_{\varepsilon,M}$
on (the Borel $\sigma$-algebra of) $\bar{\Sigma}_{M}$.

Our first result is that, no matter how large the noise level $\varepsilon\in[0,\frac{1}{2})$
is, with large enough block size $M$ it is possible to reconstruct
the starting signal $\sigma_{0}$ after performing the $M$-block
self-correction at each step.

\bigskip\noindent \textbf{Theorem} \textbf{2.1 $\forall\varepsilon\in[0,\frac{1}{2})$}
$\exists\bar{M}:\forall M>\bar{M}$ \[
\liminf_{n}\Delta_{n}(\bar{P}_{\varepsilon,M})>0.\]

\bigskip\noindent\textbf{Proof.} We first calculate the error rate
$\bar{\varepsilon}_{M}$ on the renormalized tree $\bar{T}^{(r)}$:
let $B$ be any block of size $M$ of direct descendant of some site
$v'\in B'$, where $B$ is a descendant of $B'$ in $\bar{T}^{(r)}$;
then

\begin{equation}
\bar{\varepsilon}_{M}=P_{\varepsilon}(\sum_{v\in B}\sigma_{v}<0|\sigma_{v'}=1)+\frac{1}{2}P_{\varepsilon}(\sum_{v\in B}\sigma_{v}=0|\sigma_{v'}=1).\label{eq:2.4}\end{equation}
Given $\sigma_{v'}$, the $\sigma_{v}$'s are $\{-1,1\}$-i.i.d. random
variables with $P_{\varepsilon}(\sigma_{v}=1|\sigma_{v'}=1)=1-\varepsilon>\frac{1}{2}$,
so that by large deviations theory there exists $c_{\varepsilon}>0$
such that $\bar{\varepsilon}_{M}\leq e^{-c_{\varepsilon}M}$ for all
$M>0$. Therefore, for $M$ large enough,

\begin{equation}
(1-2\bar{\varepsilon}_{M})^{2}r\geq(1-2e^{-c_{\varepsilon}M})^{2}r>1.\label{eq:2.5}\end{equation}
This implies that $\bar{\varepsilon}_{M}<\varepsilon_{c}$ and there
is reconstruction on the renormalized tree $\bar{T}^{(r)}$. By \cite{key-1}
this implies that for such $M$'s:

\begin{equation}
\liminf_{n}\left(\bar{P}_{\varepsilon,M}(\bar{\sigma}_{0}=1|\sum_{\bar{v}\in\bar{T}_{n}^{(r)}}\bar{\sigma}_{\bar{v}}>0)-\bar{P}_{\varepsilon,M}(\bar{\sigma}_{0}=-1|\sum_{\bar{v}\in\bar{T}_{n}^{(r)}}\bar{\sigma}_{\bar{v}}>0)\right)>0.\label{eq:new6}\end{equation}
Now, $\bar{\sigma}_{0}=1$ if $\sum_{v\in T_{\tilde{n}}^{(r)}}\sigma_{v}>0$
or, with probability $\frac{1}{2}$, if $\sum_{v\in T_{\tilde{n}}^{(r)}}\sigma_{v}=0$.
Therefore,

\begin{equation}
\liminf_{n}\left(P_{\varepsilon}(\sum_{v\in T_{\tilde{n}}^{(r)}}\sigma_{v}>0|\sum_{\bar{v}\in\bar{T}_{n}^{(r)}}\bar{\sigma}_{\bar{v}}>0)-P_{\varepsilon}(\sum_{v\in T_{\tilde{n}}^{(r)}}\sigma_{v}<0|\sum_{\bar{v}\in\bar{T}_{n}^{(r)}}\bar{\sigma}_{\bar{v}}>0)\right)>0.\label{eq:2.1}\end{equation}

We now show that by reading the block variables $\bar{\sigma}_{\bar{v}}$
for $\bar{v}\in\bar{T}_{n}^{(r)}$ one can reconstruct $\sigma_{0}$.
To this purpose let

\begin{eqnarray*}
A & = & \{\sigma_{0}=+1\},\\
B & = & \{\sum_{v\in T_{\tilde{n}}^{(r)}}\sigma_{v}>0\}\end{eqnarray*}

and

\begin{eqnarray*}
C & = & \{\sum_{\bar{v}\in\bar{T}_{n}^{(r)}}\bar{\sigma}_{\bar{v}}>0\}.\end{eqnarray*}

\begin{equation}
\end{equation}
We then have, by total probabilities theorem, the Markov property
and the fact that $P(A|B^{c})=P(A^{c}|B)$ (with the same equality
when $A$ and $A^{c}$ are exchanged),

\begin{gather}
P_{\varepsilon}(A|C)-P_{\varepsilon}(A^{c}|C)\nonumber \\
=P_{\varepsilon}(A|C\cap B)P_{\varepsilon}(B|C)+P_{\varepsilon}(A|C\cap B^{c})P_{\varepsilon}(B^{c}|C)\nonumber \\
-(P_{\varepsilon}(A^{c}|C\cap B)P_{\varepsilon}(B|C)+P_{\varepsilon}(A^{c}|C\cap B^{c})P_{\varepsilon}(B^{c}|C))\\
=P_{\varepsilon}(A|B)P_{\varepsilon}(B|C)+P_{\varepsilon}(A|B^{c})P_{\varepsilon}(B^{c}|C)\nonumber \\
-(P_{\varepsilon}(A^{c}|B)P_{\varepsilon}(B|C)+P_{\varepsilon}(A^{c}|B^{c})P_{\varepsilon}(B^{c}|C))\nonumber \\
=(P_{\varepsilon}(A|B)-P_{\varepsilon}(A^{c}|B))(P_{\varepsilon}(B|C)-P_{\varepsilon}(B^{c}|C))>0;\nonumber \end{gather}
the last inequality holds since it follows from (\ref{eq:2.1}) that
if $M$ is large enough, $\liminf_{n}(P_{\varepsilon}(B|C)-P_{\varepsilon}(B^{c}|C))>0$,
and it follows from the next Lemma that $P_{\varepsilon}(A|B)-P_{\varepsilon}(A^{c}|B)>0$
for every $\tilde{n}$. $\blacksquare$

\bigskip\noindent \textbf{Lemma 2.2} \emph{Consider any tree} $T^{(r)}$
\emph{and a transmission problem described by the distribution} $P_{\varepsilon}$\emph{,
let} $S_{n}(\sigma)=S_{n}=\sum_{v\in T_{n}^{(r)}}\sigma_{v}$\emph{.
Then} 

\begin{description}
\item [i)]\[
P_{\varepsilon}(S_{n-1}>0|S_{n}>0)-P_{\varepsilon}(S_{n-1}<0|S_{n}>0)>0\]

\item [ii)]\[
P_{\varepsilon}(S_{n}>0|D)-P_{\varepsilon}(S_{n}<0|D)>0\]

\end{description}
\emph{for every} $D\subseteq\{-1,1\}^{T_{n-1}^{(r)}}$ \emph{such
that} $\forall\sigma\in D$\emph{,} $S_{n-1}(\sigma)>0$.

\begin{description}
\item [iii)]\[
P_{\varepsilon}(S_{n}>0|\hat{\sigma}_{T_{n-1}^{(r)}})-P_{\varepsilon}(S_{n}<0|\hat{\sigma}_{T_{n-1}^{(r)}})>0\]

\end{description}
\emph{for every configuration} $\hat{\sigma}_{T_{n-1}^{(r)}}\in\{-1,1\}^{T_{n-1}^{(r)}}$
\emph{such that} $\sum_{v\in T_{n-1}^{(r)}}\hat{\sigma}_{v}=l>0$.

\begin{description}
\item [iv)]\[
P_{\varepsilon}(S_{n-k}>0|S_{n}>0)-P_{\varepsilon}(S_{n-k}<0|S_{n}>0)>0\]

\end{description}
\emph{for every $k=1,...,n$.}

\bigskip\noindent\textbf{Proof.} Clearly ii) implies i) taking $D=\{ S_{n-1}>0\}$,
and iii) implies ii) since 

\begin{multline*}
P_{\varepsilon}(S_{n}>0|S_{n-1}>0)\\
=\sum_{\hat{\sigma}_{T_{n-1}^{(r)}}:\sum_{v\in T_{n-1}^{(r)}}\hat{\sigma}_{v}>0}P_{\varepsilon}(S_{n}>0|\hat{\sigma}_{T_{n-1}^{(r)}})P_{\varepsilon}(\hat{\sigma}_{T_{n-1}^{(r)}}|S_{n-1}>0)\end{multline*}

\begin{equation}
\end{equation}
To show iii) assume $\sum_{v\in T_{n-1}^{(r)}}\hat{\sigma}_{v}=l>0$.
Then $S_{n}=\sum_{i=1}^{\frac{r^{n-1}-l}{2}}X_{i}+\sum_{i=\frac{r^{n-1}-l}{2}}^{\frac{r^{n-1}-l}{2}}Y_{i}+\sum_{i=r^{n-1}-l+1}^{r^{n-1}}X_{i}$
with $X_{i}$ i.i.d, $Y_{i}$ i.i.d, $X_{i},Y_{i}\in\{-r,r\}$ and
$X_{i}=\sum_{j=1}^{r}\tilde{X}_{j}$, $\tilde{X}_{j}$ i.i.d, $\tilde{X}_{j}\in\{-1,1\}$,
$P(\tilde{X}_{j}=1)=1-\varepsilon$ and $Y_{i}=\sum_{j=1}^{r}\tilde{Y}_{j}$,
$\tilde{Y}_{j}$ i.i.d, $Y_{j}\in\{-1,1\}$, $P(\tilde{Y}_{j}=1)=\varepsilon$,
all these variables being independent. So $X_{i}$ is distributed
like $S_{1}$ conditioned to $\sigma_{0}=1$ and, by symmetry of the
distribution of $S_{1}$, $X_{i}=^{d}-Y_{i}$, so that

\begin{equation}
\bar{S}_{n}=\sum_{i=1}^{\frac{r^{n-1}-l}{2}}X_{i}+\sum_{i=\frac{r^{n-1}-l}{2}+1}^{r^{n-1}-l}Y_{i}\end{equation}
is a symmetric random variable. Therefore,

\begin{multline*}
P_{\varepsilon}(S_{n}>0|\hat{\sigma}_{T_{n-1}^{(r)}})=P_{\varepsilon}(\bar{S}_{n}+\sum_{i=r^{n-1}-l+1}^{r^{n-1}}X_{i}>0)\\
=\sum_{l_{1}=l-r^{n-1}}^{r^{n-1}-l}P_{\varepsilon}(\bar{S}_{n}+\sum_{i=r^{n-1}-l+1}^{r^{n-1}}X_{i}>0|\bar{S}_{n}=l_{1})P_{\varepsilon}(\bar{S}_{n}=l_{1})\\
=\sum_{l_{1}>0}^{r^{n-1}-l}\left[P_{\varepsilon}(\sum_{i=r^{n-1}-l+1}^{r^{n-1}}X_{i}>-l_{1}|\bar{S}_{n}=l_{1})\right.\\
\left.+P_{\varepsilon}(\sum_{i=r^{n-1}-l+1}^{r^{n-1}}X_{i}>l_{1}|\bar{S}_{n}=-l_{1})\right]P_{\varepsilon}(\bar{S}_{n}=l_{1})\\
+P_{\varepsilon}(\sum_{i=r^{n-1}-l+1}^{r^{n-1}}X_{i}>0|\bar{S}_{n}=0)P_{\varepsilon}(\bar{S}_{n}=0)\\
=\sum_{l_{1}>0}^{r^{n-1}-l}\left[P_{\varepsilon}(\sum_{i=1}^{l}X_{i}>-l_{1})+P_{\varepsilon}(\sum_{i=1}^{l}X_{i}>l_{1})\right]P_{\varepsilon}(\bar{S}_{n}=l_{1})\\
+P_{\varepsilon}(\sum_{i=1}^{l}X_{i}>0)P_{\varepsilon}(\bar{S}_{n}=0)\end{multline*}

\begin{equation}
\end{equation}
By the analogous expression for $S_{n}<0$ we then need

\begin{multline}
P_{\varepsilon}(\sum_{i=1}^{l}X_{i}>-l_{1})+P_{\varepsilon}(\sum_{i=1}^{l}X_{i}>l_{1})\\
>P_{\varepsilon}(\sum_{i=1}^{l}X_{i}<-l_{1})+P_{\varepsilon}(\sum_{i=1}^{l}X_{i}<l_{1})\label{eq:12}\end{multline}
For every $l\geq1$ and $l_{1}\geq0$, we have $\sum_{i=1}^{l}X_{i}=\sum_{j=1}^{rl}\tilde{X}_{j}$
and 

\begin{equation}
P_{\varepsilon}(\sum_{j=1}^{rl}\tilde{X}_{j}>l_{1})=\sum_{h=\frac{rl+l_{1}}{2}}^{rl}\left(\begin{array}{c}
rl\\
h\end{array}\right)(1-\varepsilon)^{h}\varepsilon^{rl-h}\end{equation}
Also, by the change of variable $rl-h'=h$,

\begin{eqnarray*}
P_{\varepsilon}(\sum_{j'=1}^{rl}\tilde{X}_{j'}<-l_{1}) & = & \sum_{h'=0}^{\frac{rl-l_{1}}{2}}\left(\begin{array}{c}
rl\\
h'\end{array}\right)(1-\varepsilon)^{h'}\varepsilon^{rl-h'}\\
 & = & \sum_{h=\frac{rl+l_{1}}{2}}^{rl}\left(\begin{array}{c}
rl\\
h\end{array}\right)(1-\varepsilon)^{rl-h}\varepsilon^{h}\end{eqnarray*}

\begin{equation}
\end{equation}
So that, for $\varepsilon<\frac{1}{2}$,

\begin{multline}
P_{\varepsilon}(\sum_{i=1}^{l}X_{i}>l_{1})-P_{\varepsilon}(\sum_{i=1}^{l}X_{i}<-l_{1})\\
\sum_{h=\frac{rl+l_{1}}{2}}^{rl}\left(\begin{array}{c}
rl\\
h\end{array}\right)(1-\varepsilon)^{rl-h}\varepsilon^{rl-h}((1-\varepsilon)^{2h-rl}-\varepsilon^{2h-rl})>0.\\
\label{eq:13}\end{multline}
This shows (\ref{eq:12}) since we have seen one strict inequality
between two terms, and the other two terms satisfy

\begin{multline}
P_{\varepsilon}(\sum_{i=1}^{l}X_{i}>-l'_{1})-P_{\varepsilon}(\sum_{i=1}^{l}X_{i}<l'_{1})\\
=P_{\varepsilon}(\sum_{i=1}^{l}X_{i}\geq l'_{1})-P_{\varepsilon}(\sum_{i=1}^{l}X_{i}\leq-l'_{1})>0\\
\end{multline}
for the same inequality (\ref{eq:13}) applied to $l_{1}=l'_{1}-1\geq0$.

Finally, (iv) is shown using iteratively (\ref{eq:new6}) for $k$
larger than one with

\begin{eqnarray*}
A & = & \{ S_{n-k}>0\}\\
B & = & \{ S_{n-k+1}>0\}\end{eqnarray*}
and

\begin{eqnarray*}
C & = & \{ S_{n}\geq0\}.\end{eqnarray*}
$\blacksquare$

\section{Within-descent self-correction: strict inequality of critical points}

Our aim is to consider within-descent self-correction at some level
$k$. To this purpose we take a vertex $v$ in some generation $mk$,
$m\in\mathbb{N}$, and look at its $r^{k}$ descendants $k$ generation
down (thus in $T_{(m+1)k}^{(r)}$) as generated by the transmission;
we then force all such descendants to agree to their majority (with
random choice if there is no majority). Transmission is then resumed
as usual from the modified status. This amounts to define a map $\Phi_{k}:\{-1,1\}^{T^{(r)}}\rightarrow\Sigma_{k}$
given by

\[
\Phi_{k}(\sigma)_{v}=\left\{ \begin{array}{c}
\begin{array}{c}
\:\;1\quad\textrm{with probability }1\quad\textrm{if}\quad\sum_{s'_{2}=1}^{r^{k}}\sigma_{mk,s_{1}r^{k}+s'_{2}}>0\\
\begin{array}{c}
\\\\\end{array}\end{array}\\
\begin{array}{c}
-1\quad\textrm{with probability }1\quad\textrm{if}\quad\sum_{s'_{2}=1}^{r^{k}}\sigma_{mk,s_{1}r^{k}+s'_{2}}<0\\
\begin{array}{c}
\\\\\end{array}\end{array}\\
\quad\left\{ \begin{array}{c}
\:\;1\quad\textrm{with probability }1/2\\
-1\quad\textrm{with probability }1/2\end{array}\quad\right.\textrm{if}\quad\sum_{s'_{2}=1}^{r^{k}}\sigma_{mk,s_{1}r^{k}+s'_{2}}=0\end{array}\right.\]
\begin{equation}
\end{equation}
if $v\in T_{mk}^{(r)}$, with $v=(mk,s_{1}r^{k}+s_{2})$, $s_{1}=0,...,r^{k(m-1)}-1$,
$s_{2}=1,...,r^{k}$; otherwise

\begin{equation}
\Phi_{k}(\sigma)_{v}=\sigma_{v}.\label{eq:}\end{equation}
 As before, the transmission is self-corrected by the map $\Phi_{k}$
every $k$ steps: $\sigma_{T_{\rightarrow mk}^{(r)}}\in\Sigma_{k}$
generates $\sigma_{T_{\rightarrow(m+1)k}^{(r)}}\in\{-1,1\}^{T^{(r)}}$
as usual, and then we take $\Phi_{k}(\sigma_{T_{\rightarrow(m+1)k}^{(r)}})\in\Sigma_{k}$.
The distribution $P_{\varepsilon}^{(k)}$ of the self-corrected configuration
is then recursively defined by \begin{equation}
P_{\varepsilon}^{(k)}(\sigma_{T_{\rightarrow(m+1)k}^{(r)}}|\sigma_{T_{\rightarrow mk}^{(r)}})=P_{\varepsilon}(\Phi_{k}^{-1}\sigma_{T_{\rightarrow(m+1)k}^{(r)}\setminus T_{\rightarrow mr}^{(r)}}|\sigma_{T_{\rightarrow mk}^{(r)}}).\end{equation}
Notice that $P_{\varepsilon}^{(k)}$ is no longer a Markov chain but
the conditional probabilities satisfy

\begin{eqnarray*}
P_{\varepsilon}^{(k)}(\sigma_{T_{n}^{(r)}}|\sigma_{T_{\rightarrow(n-1)}^{(r)}}) & = & P_{\varepsilon}(\sigma_{T_{n}^{(r)}}|\sigma_{T_{\rightarrow(n-1)}^{(r)}})\\
 & = & P_{\varepsilon}(\sigma_{T_{n}^{(r)}}|\sigma_{T_{(n-1)}^{(r)}})\end{eqnarray*}

\begin{equation}
\label{eq:new25}\end{equation}
for all $n$ not of the form $n=mk$.

Next, for $\sigma\in\Sigma_{k}$, let $\Psi_{k}(\sigma)\in T^{(r^{k})}$
be defined by

\begin{equation}
\Psi_{k}(\sigma)_{v}=\sigma_{(mk,s_{1}r^{k}+1)}\label{eq:5}\end{equation}
if $v\in T_{m}^{(r^{k})}$, $v=(m,s_{1}r^{k}+s_{2})$, $s_{1}=0,...,r^{k(m-1)}-1$,
$s_{2}=1,...,r^{k}$. Note that $\Psi_{k}(\sigma)$ is a configuration
of an almost regular tree $T^{(r^{k})}$: $T^{(r^{k})}$ has branching
rate $1$ at the starting vertex and then $r^{k}$ at all other vertices.
As we will see, the initial segment makes no difference in our arguments,
and, therefore, we adopt the slight abuse of notation $T^{(r^{k})}$
(which in our definitions indicates a regular tree).

Using $P_{\varepsilon}^{(k)}$ we define the self-corrected critical
distortions 

\begin{equation}
\varepsilon_{c,r}(k)=\sup\{\varepsilon:\liminf_{n}\Delta_{n}(P_{\varepsilon}^{(k)})>0\}.\end{equation}

Note that on $\Psi_{k}(\Phi_{k}(\{-1,1\}^{T^{(r)}}))=T^{(r^{k})}$
the distribution $\Psi_{k}(P_{\varepsilon}^{(k)})=P_{\varepsilon}^{(k)}\cdot\Psi_{k}^{-1}$
is a Markov chain, by the definition of $P_{\varepsilon}^{(k)}$,
and thus it is again a transmission model with error rate $\varepsilon(k)$.
In other words, $\Psi_{k}(P_{\varepsilon}^{(k)})=P_{\varepsilon(k)}$
on $T^{(r^{k})}$.

We first show that reconstruction under $\Psi_{k}(P_{\varepsilon}^{(k)})$
on $T^{(r^{k})}$ is equivalent to reconstruction under the $k$-self
corrected distribution $P_{\varepsilon}^{(k)}$.

\bigskip\noindent\textbf{Lemma 3.1} $\liminf_{n}\Delta_{n}(P_{\varepsilon}^{(k)})>0$
\emph{if and only if} $\liminf_{n}\Delta_{n}(\Psi_{k}(P_{\varepsilon}^{(k)}))>0$ 

\bigskip\noindent\textbf{Proof.} First, observe that $\liminf_{n}\Delta_{n}(\Psi_{k}(P_{\varepsilon}^{(k)}))>0$
on $\Psi_{k}(\Sigma_{k})$ if and only if $\liminf_{n}\Delta_{n}(\Psi_{k}(P_{\varepsilon}^{(k)}))>0$
on $T^{(r^{k})}$. In fact, on $\Psi_{k}(\Sigma_{k})$ we obtain

\begin{gather}
P_{\varepsilon(k)}(S_{n}>0|\sigma_{0}>0)\nonumber \\
=\liminf_{n}\left[(1-\varepsilon(k))P_{\varepsilon(k)}(S_{n}>0|\sigma_{(1,1)}>0)\right.\\
\left.+\varepsilon(k)P_{\varepsilon(k)}(S_{n}>0|\sigma_{(1,1)}<0)\right]\nonumber \\
=(1-2\varepsilon(k))\liminf_{n}P_{\varepsilon(k)}(S_{n}>0|\sigma_{(1,1)}>0)+\varepsilon(k)\nonumber \\
\end{gather}
so that 

\begin{gather*}
P_{\varepsilon(k)}(S_{n}>0|\sigma_{0}>0)-P_{\varepsilon(k)}(S_{n}<0|\sigma_{0}>0)\\
=(1-2\varepsilon(k))\Delta_{n}(\Psi_{k}(P_{\varepsilon}^{(k)}));\end{gather*}

\begin{equation}
\end{equation}
the $\liminf_{n}$ of the last expression is positive if and only
if $\liminf_{n}\Delta_{n}(\Psi_{k}(P_{\varepsilon}^{(k)}))>0$ on
$T^{(r^{k})}$ as $\varepsilon(k)<1/2$. Now, observe that $\liminf_{n}\Delta_{n}(P_{\varepsilon}^{(k)})>0$
implies $\liminf_{mk}\Delta_{mk}(P_{\varepsilon}^{(k)})>0$, that
is $\liminf_{n}\Delta_{n}(\Psi_{k}(P_{\varepsilon}^{(k)}))>0$ on
$\Psi_{k}(\{-1,1\}^{T(r)})$. 

To show the reverse implication, notice that for every level $n$
of $T^{(r^{k})}$ not of the form $n=mk$ we have

\begin{gather}
\Delta_{n}(P_{\varepsilon}^{(k)})=\Delta_{n-1}(P_{\varepsilon}^{(k)})(P_{\varepsilon}^{(k)}(S_{n}>0|S_{n-1}>0)-P_{\varepsilon}^{(k)}(S_{n}<0|S_{n-1}>0))\end{gather}
where if $n-1=mk$ then $S_{n-1}=\sum_{v\in T_{n-1}^{(r)}}(\Phi_{k}(\sigma))_{v}$.
In all cases, the event $D=\{ S_{n-1}>0\}$ is such that $\hat{\sigma}\in D$
satisfies $\sum_{v\in T_{n-1}^{(r)}}\hat{\sigma}_{v}>0$; this implies
$P_{\varepsilon}^{(k)}(S_{n}>0|S_{n-1}>0)>P_{\varepsilon}^{(k)}(S_{n}<0|S_{n-1}>0)$
by part ii) of Lemma 2.2 applied to $P_{\varepsilon}^{(k)}$, since,
by (\ref{eq:new25}), the conditional probabilities coincide with
those of $P_{\varepsilon}$. 

Therefore, computing $\Delta_{n}(P_{\varepsilon}^{(k)})$ by finite
iteration from the maximum level $mk<n$, $\liminf_{mk}\Delta_{mk}(P_{\varepsilon}^{(k)})>0$
implies $\liminf_{n}\Delta_{n}(P_{\varepsilon}^{(k)})>0$. $\blacksquare$

Our next aim is to show that $\varepsilon_{c,r}(k)>\bar{\varepsilon}_{c,r}$,
which is to say $p_{c,r}(k)<\bar{p}_{c,r}$, where $\bar{\varepsilon}_{c,r}$
is the critical distortion rate for majority or maximum likelihood
reconstruction on $T^{(r)}$.

In order to do this we introduce another random transformation, the
fraction identification transform $\tilde{\Phi}_{k}:\{-1,1\}^{T^{(r)}}\rightarrow\{-1,1\}^{T^{(r)}}$
given by

\begin{equation}
\tilde{\Phi}_{k}(\sigma)_{v}=\sigma_{\bar{v}}\label{eq:}\end{equation}
if $v\in T_{mk}^{(r)}$, with $v=(mk,s_{1}r^{k}+s_{2})$, $s_{1}=0,...,r^{k(m-1)}-1$,
$s_{2}=1,...,r^{k}$, and $\bar{v}=(mk,s_{1}r^{k}+\bar{s}_{2})$,
$\bar{s}_{2}=1,...,r^{k}$ uniformly chosen at random. Otherwise \begin{equation}
\tilde{\Phi}_{k}(\sigma)_{v}=\sigma_{v}.\label{eq:}\end{equation}

As before, for $\sigma\in\tilde{\Phi}_{k}(\{-1,1\}^{T^{(r)}})$, let
$\tilde{\Psi}_{k}(\sigma)\in T^{(r^{k})}$ be defined by

\begin{equation}
\tilde{\Psi}_{k}(\sigma)_{v}=\sigma_{(m,s_{1}r^{k})}.\label{eq:6}\end{equation}

Now, the strict inequality between the self-corrected critical distortion
and the original one can be proven. The strict inequality holds for
all values of $k$ and $r$ except for the one step correction on
binary trees. 

\bigskip\noindent\textbf{Theorem 3.2} \emph{If} $k>1$ \emph{or}
$k=1$\emph{,} $r>2$

\begin{eqnarray}
\varepsilon_{c,r}(k) & > & \bar{\varepsilon}_{c,r};\label{eq:}\end{eqnarray}
\begin{eqnarray}
\varepsilon_{c,2}(1) & = & \bar{\varepsilon}_{c,2}\label{eq:}\end{eqnarray}

\bigskip         To prove this fact, we explicitly compute the noise
change under the fraction identification. On $\tilde{\Phi}_{k}(\{-1,1\}^{T(r)})$
the probability distribution $\tilde{P}_{\varepsilon}(k)$ which implements
the fraction transform is defined as $P_{\varepsilon}^{(k)}$ with
$\varepsilon(k)$ replaced by $\tilde{\varepsilon}(k)=1-\frac{1}{r^{k}}\sum_{s'_{2}=1}^{r^{k}}(2\sigma_{mk,s_{1}r^{k}+s'_{2}}-1)$.
Note that $\Psi_{k}(\tilde{P_{\varepsilon}}^{(k)})=P_{\tilde{\varepsilon}(k)}$
on $T^{(r^{k})}$. We then have

\bigskip\noindent\textbf{Lemma 3.3} $\forall\varepsilon$, $\forall k$

\begin{equation}
1-2\tilde{\varepsilon}(k)=(1-2\varepsilon)^{k}\label{eq:1}\end{equation}
 \emph{therefore the critical distortion} $\tilde{\varepsilon}_{c,r}(k)=\sup\{\varepsilon:\liminf_{n}\Delta_{n}(\tilde{P}_{\varepsilon}^{(k)})>0\}$
\emph{equals} $\varepsilon_{c,r^{k}}$.

\bigskip\noindent\textbf{Proof.} Denote by $X_{k}$ the number of
$1$'s at level $k$. By definition and linearity of expected values,\begin{equation}
\tilde{\varepsilon}_{k}(k)=1-\frac{1}{r^{k}}E_{\varepsilon}(X_{k}|\sigma_{0}=1)=1-P_{\varepsilon}(\sigma_{\bar{v}}=1|\sigma_{0}=1)\label{eq:}\end{equation}
for every $\bar{v}\in T_{k}^{(r)}$. The last probability refers to
a one-dimensional Markov chain of length $k$ with distortion probability
$\varepsilon$, and can be easily computed. Alternatively, (\ref{eq:1})
can be verified by induction, since by the last equality, $1-2\tilde{\varepsilon}(1)=1-2\varepsilon$
and \begin{equation}
\tilde{\varepsilon}(k)=\varepsilon(1-\tilde{\varepsilon}(k-1))+(1-\varepsilon)\tilde{\varepsilon}(k-1),\label{eq:}\end{equation}
so that \begin{equation}
1-2\tilde{\varepsilon}(k)=(1-2\varepsilon)(1-2\tilde{\varepsilon}(k-1))=(1-2\varepsilon)^{k}.\label{eq:}\end{equation}
From \cite{key-1}, $(1-2\varepsilon_{c,r})^{2}r=1$ and since $\Psi_{k}(\tilde{P}_{\varepsilon}(k))$
is on $T^{(r^{k})}$, on this second tree criticality is identified
by $(1-2\varepsilon_{c,r^{k}})^{2}r^{k}=1$ and (\ref{eq:1}) implies
$(1-2\tilde{\varepsilon}_{c,r}(k))^{2}r^{k}=((1-2\varepsilon_{c,r})^{k})^{2}r^{k}=((1-2\varepsilon_{c,r})^{2}r)^{k}=1$
. So $\tilde{\varepsilon}_{c,r}(k)=\varepsilon_{c,r^{k}}$. $\blacksquare$

\bigskip\noindent\textbf{Proof of Theorem 3.2} Introduce 

\begin{equation}
T_{k,r}(\varepsilon)=\frac{1}{r^{k}}\sum_{l=0}^{\frac{r^{k}-1}{2}}l\left(P_{\varepsilon}(X_{k}=l|\sigma_{0}=0)-P_{\varepsilon}(X_{k}=l|\sigma_{0}=1)\right)\label{eq:}\end{equation}
when $r$ is odd, and

\begin{equation}
T_{k,r}(\varepsilon)=\frac{1}{r^{k}}\sum_{l=0}^{\frac{r^{k}}{2}-1}l\left(P_{\varepsilon}(X_{k}=l|\sigma_{0}=0)-P_{\varepsilon}(X_{k}=l|\sigma_{0}=1)\right)\label{eq:}\end{equation}
when $r$ is even\emph{.} For $r$ odd, we have

\begin{eqnarray*}
T_{k,r}(\varepsilon) & = & \frac{1}{r^{k}}\sum_{l=\frac{r^{k}+1}{2}}^{r^{k}}(r^{k}-l)P_{\varepsilon}(X_{k}=r^{k}-l|\sigma_{0}=0)-\frac{1}{r^{k}}\sum_{l=0}^{\frac{r^{k}-1}{2}}lP_{\varepsilon}(X_{k}=l|\sigma_{0}=1)\\
 & = & \tilde{\varepsilon}(k)-\varepsilon(k)\end{eqnarray*}
\begin{equation}
\end{equation}
and, for $r$ even 

\begin{eqnarray*}
T_{k,r}(\varepsilon) & = & \frac{1}{r^{k}}\sum_{l=\frac{r^{k}}{2}+1}^{r^{k}}(r^{k}-l)P_{\varepsilon}(X_{k}=r^{k}-l|\sigma_{0}=0)\\
 &  & -\frac{1}{r^{k}}\sum_{l=0}^{\frac{r^{k}-1}{2}}lP_{\varepsilon}(X_{k}=l|\sigma_{0}=1)+\frac{1}{2}P_{\varepsilon}(X_{k}=\frac{r^{k}}{2})-\frac{1}{2}P_{\varepsilon}(X_{k}=\frac{r^{k}}{2})\\
 & = & \tilde{\varepsilon}(k)-\varepsilon(k)\end{eqnarray*}

\begin{equation}
\end{equation}
By Lemma 3.3 $\tilde{\varepsilon}_{c,r}(k)=\bar{\varepsilon}_{c,r^{k}}$
and $T_{1,2}(\bar{\varepsilon}_{c,2})=0$, so it is sufficient to
show that $T_{k,r}(\bar{\varepsilon}_{c,r})>0$ for the non trivial
cases of $k$ and $r$. Theorem 1.4 in \cite{key-1} shows that $P_{\varepsilon}(X_{k}=l|\sigma_{0}=0)\geq P_{\varepsilon}(X_{k}=l|\sigma_{0}=1)$
if $r^{k}-l>l$. To have strict inequality it is sufficient to show
that $P_{\varepsilon}(X_{k}=1|\sigma_{0}=0)>P_{\varepsilon}(X_{k}=1|\sigma_{0}=1)$.
This will be done by induction in $k$. We focus on the number $i$
of distortions of $\sigma_{0}$ at the first step. The index $i$
runs from $0$ to $r$, but it is convenient to group together the
$i$-th and the $(r-i)$-th terms. Note that $P_{\varepsilon}(X_{1}=i|\sigma_{0}=0)=\left(\begin{array}{c}
r\\
i\end{array}\right)\varepsilon^{i}(1-\varepsilon)^{r-i}$. Assuming $\bar{i}=\frac{r+1}{2}$ for $r$ odd and $\bar{i}=\frac{r}{2}+1$
if $r$ is even and $i\geq\bar{i}$, the terms in $T_{k,r}$ can be
collected like this

\begin{eqnarray}
T_{k,r}(\varepsilon) & = & \sum_{i=\bar{i}}^{r}\left(\begin{array}{c}
r\\
i\end{array}\right)T_{k,r,i}(\varepsilon)\label{eq:}\end{eqnarray}
with \begin{eqnarray*}
T_{k,r,i}(\varepsilon) & = & \left[\varepsilon^{i}(1-\varepsilon)^{r-i}-(1-\varepsilon)^{i}\varepsilon^{r-i}\right]\\
 &  & \cdot\left[iP_{\varepsilon}(X_{k-1}=1|\sigma_{0}=1)(P_{\varepsilon}(X_{k-1}=0|\sigma_{0}=1))^{i-1}\right.\\
 &  & \cdot(P_{\varepsilon}(X_{k-1}=0|\sigma_{0}=0))^{r-i}+(r-i)(P_{\varepsilon}(X_{k-1}=0|\sigma_{0}=1))^{i}\\
 &  & \cdot P_{\varepsilon}(X_{k-1}=1|\sigma_{0}=0)(P_{\varepsilon}(X_{k-1}=0|\sigma_{0}=0))^{r-i-1}\\
 &  & -iP_{\varepsilon}(X_{k-1}=1|\sigma_{0}=0)(P_{\varepsilon}(X_{k-1}=0|\sigma_{0}=0))^{i-1}\\
 &  & \cdot(P_{\varepsilon}(X_{k-1}=0|\sigma_{0}=1))^{r-i}-(r-i)(P_{\varepsilon}(X_{k-1}=0|\sigma_{0}=0))^{i}\\
 &  & \left.\cdot P_{\varepsilon}(X_{k-1}=1|\sigma_{0}=1)(P_{\varepsilon}(X_{k-1}=0|\sigma_{0}=1))^{r-i-1}\right]\end{eqnarray*}

\begin{equation}
\end{equation}
Now, the first factor is negative if $\varepsilon\in(0,1/2)$ in particular
if $\varepsilon=\bar{\varepsilon}_{c,r}$. We now show that the second
factor is negative as well under the hypothesis that the statement
is true for $k-1$. 

The $(r-i)$ terms of the second addend are greater than or equal
to $(r-i)\leq i$ terms taken from the third addend since \begin{equation}
P_{\varepsilon}(X_{k-1}=0|\sigma_{0}=0)\geq P_{\varepsilon}(X_{k-1}=0|\sigma_{0}=1)\label{eq:}\end{equation}
 again by \cite{key-1}. The remaining $(2i-r)$ terms from the third
addend are strictly less than $(2i-r)\leq i$ terms taken from the
first since 

\begin{multline}
P_{\varepsilon}(X_{k-1}=1|\sigma_{0}=0)P_{\varepsilon}(X_{k-1}=0|\sigma_{0}=0)\\
>P_{\varepsilon}(X_{k-1}=1|\sigma_{0}=1)P_{\varepsilon}(X_{k-1}=0|\sigma_{0}=1);\\
\end{multline}
in fact, $P_{\varepsilon}(X_{k-1}=0|\sigma_{0}=0)\geq P_{\varepsilon}(X_{k-1}=0|\sigma_{0}=1)$
follows from \cite{key-1}, and $P_{\varepsilon}(X_{k-1}=1|\sigma_{0}=0)>P_{\varepsilon}(X_{k-1}=1|\sigma_{0}=1)$
follows by the induction hypothesis. 

Finally, the remaining $(r-i)$ terms in the first addend are greater
than or equal to the $(r-i)$ terms in the fourth addend again by
\cite{key-1}.

For $r=2$ and $k=2$ the statement is true, as, by direct computation,
we have, for some $f(\varepsilon)$,

\begin{multline}
P_{\varepsilon}(X_{2}=1|\sigma_{0}=0)-P_{\varepsilon}(X_{2}=1|\sigma_{0}=1)\\
=4(1-\varepsilon)^{5}\varepsilon+2\varepsilon(1-\varepsilon)f(\varepsilon)+4(1-\varepsilon)\varepsilon^{5}-8(1-\varepsilon)^{3}-2\varepsilon(1-\varepsilon)f(\varepsilon)\\
=4\varepsilon(1-\varepsilon)((1-\varepsilon)^{2}+\varepsilon^{2})^{2}\\
\end{multline}

which is positive for $\varepsilon\in(0,1/2)$. For $r>2$ and $k=1$
the statement is true as well in the same domain as $P_{\varepsilon}(X_{1}=1|\sigma_{0}=0)-P_{\varepsilon}(X_{1}=1|\sigma_{0}=1)=\varepsilon(1-\varepsilon)^{r-1}>(1-\varepsilon)\varepsilon^{r-1}$.
$\blacksquare$

\section{Limit of within-descent self-corrected critical distortions}

The transmission model we are considering can equivalently be rewritten
(see \cite{key-1}) as an Ising model $\mu_{\beta}$ with inverse
temperature $\beta$ such that \begin{equation}
\tanh(\beta)=1-2\varepsilon\end{equation}
 and \begin{equation}
\mu_{\beta,\eta}(\sigma_{T_{\rightarrow n}^{(r)}})=\frac{1}{Z}e^{-\beta\sum_{(^{\leftarrow}v,v)}\sigma_{^{\leftarrow}v}\sigma_{v}}\end{equation}
 where $\mu_{\beta}$ is any weak limit of $\mu_{\beta,\eta}$. In
turn, this can be represented as an FK model, see \cite{key-5}. The
usual FK parameter $p'=1-e^{2\beta}$ can then be modified on the
tree, to account also for the number of clusters, to $p=\frac{p'}{2-p'}=\tanh(\beta)=1-2\varepsilon$.
With $H=\{0,1\}^{\mathcal{E}(T^{(r)})}$, where $\mathcal{E}(T^{(r)})$
are the length $1$ edges of $T^{(r)}$ and $\eta\in H$, denoting
by $\mathcal{E}(T_{\rightarrow n}^{(r)})$ the edges of $T_{\rightarrow n}^{(r)}$,
we have

\begin{equation}
\nu_{p}(\eta_{\mathcal{E}(T_{\rightarrow n}^{(r)})})=\prod_{^{\leftarrow}v,v\in T_{\rightarrow n}^{(r)}}p^{\eta_{(^{\leftarrow}v,v)}}(1-p)^{1-\eta_{(^{\leftarrow}v,v)}}.\label{eq:}\end{equation}
Therefore, the FK model is simply an independent Galton-Watson branching
process with each descendant generated independently with probability
$p$. The relation between $\nu_{p}$ and $\mu_{\beta}$ is the usual
(see \cite{key-5})

\begin{equation}
\mu_{\beta}(\sigma_{T_{\rightarrow n}^{(r)}})=\sum_{\eta_{\mathcal{E}(T_{\rightarrow n}^{(r)})}\sim\sigma_{T_{\rightarrow n}^{(r)}}}\nu(\eta_{\mathcal{E}(T_{\rightarrow n}^{(r)})})\frac{1}{Cl(\eta_{\mathcal{E}(T_{\rightarrow n}^{(r)})})}\label{eq:}\end{equation}
where $\sim$ means that $\sigma$ is compatible with $\eta$, i.e.,
$\sigma_{^{\leftarrow}v}\sigma_{v}\eta_{(^{\leftarrow}v,v)}\geq0$,
and $Cl(\eta_{\mathcal{E}(T_{\rightarrow n}^{(r)})})$ equals the
number of $\sigma$'s compatible with the given $\eta$, i.e. the
number of site clusters determined by $1$-edges in $\eta$.

In this section we want to show that $\varepsilon_{c,r}(k)\rightarrow\varepsilon_{c,r}$,
i.e. $p_{c,r}(k)\rightarrow p_{c,r}$ and the main results will be

\bigskip\noindent\textbf{Theorem 4.1} There exist $c_{1},c_{2}>0$
and a function $\alpha_{k}>0$, $\lim_{k\rightarrow\infty}\alpha_{k}=0$
such that

\begin{equation}
\frac{1}{r}\vee\frac{1}{c_{1}^{\frac{1}{2k}}r}\leq p_{c}(k)\leq\frac{1+\alpha_{k}}{c_{2}^{\frac{1}{2k}}r}\end{equation}

\bigskip         so that it easily follows 

\bigskip\noindent\textbf{Corollary 4.2\[
\lim_{k\rightarrow\infty}p_{c}(k)=\frac{1}{r}.\]
}

\bigskip\noindent

The FK representation is thus a process in which each edge $e\in\mathcal{E}(T^{(r)})$
is open, i.e. $\eta_{e}=1$, independently of all other edges, with
probability $p$. The open edges are then just the (randomly selected)
error fre edges in the transmission, in the sense that, given the
configuration of the edges, the signal is generated by:

\begin{description}
\item [i)]fixing the signal $\sigma_{0}$ at the origin;
\item [ii)]having the signal transmitted error free through the open edges;
\item [iii)]having the signal chosen at random with equal probability through
the closed edges.
\end{description}
Seen globally, the set of vertices of $T^{(r)}$ falls apart into
maximal connected components connected by open edges, and such components
are called clusters. The cluster containing a vertex $v$ is indicated
by $C(v)$. Notice that $C(0)$ describes the descendants of a Galton-Watson
process with offspring distribution Bernoulli of parameters $r$ and
$p$. The configuration of FK edges can also be described by some
$\eta\in\{0,1\}^{\mathcal{E}(T^{(r)})}$.

As before, let $T_{n}^{(r)}$ be the vertices in the $n$-th generation
of the tree. The vertices of $T_{n}^{(r)}$ also fall apart into {}``clusters''
connected, via the entire tree, by open edges (these {}``clusters''
are just the intersection of the clusters of $T^{(r)}$ with $T_{n}^{(r)}$).
Given a configuration $\eta\in\{0,1\}^{\mathcal{E}(T^{(r)})}$ of
open, i.e. value $1$, FK edges, let $Z_{i}=Z_{i}(\eta)$, $i=1,...,m_{n}=m_{n}(\eta)$,
be the clusters of $T_{n}^{(r)}$ in $\eta$, $1\leq m_{n}\leq r^{n}$,
and let $z_{i}=|Z_{i}|$. 

Notice that $\Psi_{k}(\Phi_{k}(\sigma))$ is a configuration of $T^{(r^{k})}$
and that on such tree there is reconstruction if the FK density $p(k)=p_{r^{k}}$
is such that $p_{r^{k}}^{2}r^{k}>1$ (see \cite{key-1}). 

On the other hand, by our construction, $p_{r^{k}}=1-2P(S_{k}>0|\sigma_{0}=1)$,
so we need a lower bound for this expression. Such lower bound is
achieved by estimating the size of $C(0)\cap T_{k}^{(r)}$, which
is the set carrying information, and the value of $\sum_{i=1}^{m_{k}}Z'_{i}$,
where $Z'_{i}$ are independent symmetric random variables taking
values in $\{-z_{i},z_{i}\}$, i.e. distributed as the $Z_{i}$'s.
This last sum can be estimated via the normal approximation using
Berry-Essen estimates of the error. This, however, involves second
and third moments of $Z_{i}$, and we need to develop a somewhat elaborate
bound on these moments since simple ones based on the maximum size
of $Z_{i}$ are not sufficient.

Such bounds on the sums of moments of $Z_{i}$'s are determined in
Theorems 4.2 and 4.3 below, as follows. First, notice that in creating
the $k$-th generation roughly $(1-p)r^{k-1}$ vertices are isolated,
thus giving rise to the same number of $Z_{i}$'s taking values in
$\{-1,1\}$. Therefore, $\sum_{i=1}^{m_{k}}z_{i}^{2}\geq cr^{k}$
for some $c>0$ and our first two estimates show that this bound is
nearly optimal. On the other hand, the largest cluster is of size
roughly $(pr)^{k}$, so that $z_{i}^{3}\simeq(pr)^{3k}=(p^{2}r)^{k}(pr^{2})^{k}\leq(1-c)^{k}(pr^{2})^{k}$
if $p^{2}r<1$. Our last estimate shows that also this bound is nearly
optimal. Note that this estimate cannot hold if $p^{2}r\geq1$, so
that it provides no information about the reconstruction regime of
the original tree. 

We first need a large deviation result for the size of the set of
vertices $R_{n}=C(0)\cap T_{n}^{(r)}$, i.e. for the survival set
of the Galton-Watson process in the $n$- th generation. Let $P_{p}=P_{\varepsilon}$
for $p=1-2\varepsilon$.

\bigskip\noindent \textbf{Lemma 4.3} \emph{Let $\gamma=logr/log(pr)>1$
and $\gamma^{*}$ such that $1/\gamma+1/\gamma^{*}=1$ and let $W=\lim_{n\rightarrow\infty}\frac{|R_{n}|}{(pr)^{n}}$
(see \cite{key-6}). Indicating by $P$ the distribution of $W$ and
by $E$ the expected value with respect to $P$, if $pr>1$ then there
exist $M,c_{1},c_{2},c_{3}>0$ such that if $\varepsilon>0$ is such
that $(1+\varepsilon)^{\gamma^{*}}<(pr)^{1/3}$ and $l\in\mathbb{N}$
is such that $\frac{((1+\varepsilon)/2)^{\gamma}}{\gamma^{*}(\gamma\tau)^{1/(\gamma-1)}}\leq c_{1}(pr)^{1/3}$
and $(1+\varepsilon)^{l}/2>M\vee1$ with $\tau=\max_{x<pr}H(x)<\infty$
and $H(x)=x^{-\gamma}\log(B_{r}\cdot\Phi(x))$, $\Phi(s)=E(e^{sW})$
and $B_{r}$ the Bottcher's function (see \cite{key-13}), then}

\begin{equation}
P_{p}(|R_{l}|\geq(1+\varepsilon)^{l}p^{l}r^{l})\leq c_{2}e^{c_{3}(1+\varepsilon)^{\gamma^{*}l}}\end{equation}
\emph{for all $l\in\mathbb{N}$.}

\bigskip\noindent\textbf{Proof.} By large deviation properties of
$W$, there exists $M>0$ such that for all $x>M$

\begin{equation}
P(W\geq x)\leq\exp\left(\frac{x^{\gamma^{*}}}{\gamma^{*}(\gamma\tau)^{1/(\gamma-1)}}\right)\end{equation}
for all $x$. Also, there exist $c_{4},c_{5}>0$ such that

\begin{equation}
P\left(\left|\frac{|R_{n}|}{(pr)^{n}}-W\right|\geq1\right)\leq c_{4}e^{c_{5}(pr)^{n/3}},\end{equation}
for all $n$, see \cite{key-7}, Theorem 5; the conditions of that
result are easily met by considering a process with the offspring
of $R_{n}$ plus one additional offspring in each vertex. Therefore,
under the current assumptions, for some $c_{2}\geq c_{4}+1$ and all
$l\in\mathbb{N}$

\begin{eqnarray*}
P(|R_{n}|\geq(1+\varepsilon)^{l}p^{l}r^{l}) & \leq & P\left(\left|\frac{|R_{n}|}{(pr)^{n}}-W\right|\geq1\right)+P(W\geq(1+\varepsilon)^{l}/2)\\
 & \leq & c_{4}e^{c_{5}(pr)^{l/3}}+\exp\left(\frac{((1+\varepsilon)^{l}/2)^{\gamma^{*}}}{\gamma^{*}(\gamma\tau)^{1/(\gamma-1)}}\right)\\
 & \leq & c_{2}e^{c_{3}(1+\varepsilon)^{\gamma^{*}l}}\end{eqnarray*}

\begin{equation}
\end{equation}
if $c_{3}=\frac{1}{\gamma^{*}(2\gamma^{*}\gamma\tau)^{1/(\gamma-1)}}$.
$\blacksquare$

\bigskip\noindent\textbf{Theorem 4.4} $\forall p$ \emph{and} $r$
\emph{with} $P_{p}$\emph{-probability one there exists a constant}
$c_{7}=c_{7}(\eta)>0$ \emph{such that}

\begin{equation}
\sum_{i=1}^{m_{k}}z_{i}^{2}\geq c_{7}r^{k}\label{eq:2}\end{equation}
\emph{for all $k$ larger than some $\bar{k}_{7}(\eta)$.}

\bigskip\noindent\textbf{Proof.} $\sum_{i=1}^{m_{k}}z_{i}^{2}\geq\sum_{C:C\cap T_{k}^{(r)}\neq\emptyset,|C|=1}|C|^{2}=|\{ C\subseteq T_{k}^{(r)}:|C|=1\}|=:I_{k}$.
For every $b=(^{\leftarrow}v,v)$, $^{\leftarrow}v\in T_{k-1}^{(r)}$,
$\eta_{b}$ is independently chosen to be $0$ with probability $1-p$,
and in such a case $C(v)=\{ v\}$. So, by large deviations estimates
for $r^{k}$ i.i.d. binary random variables, if $c_{7}=\frac{1-p}{2}$,
$P(I_{k}\leq c_{7}r^{k})\leq e^{-c_{3}\frac{(1-p)}{2}r^{k}}$ for
some $c_{3}>0$ (see, for instance \cite{key-9}) 

Therefore, $\sum_{k=1}^{\infty}P_{p}(\eta:I_{k}\leq c_{7}r^{k})\leq\sum_{k=1}^{\infty}e^{-c_{3}\frac{(1-p)}{2}r^{k}}<\infty$
and by Borel-Cantelli the statement holds with $P_{p}$-probability
$1$ for large $k$ with $c_{7}=\frac{1-p}{2}$. $\blacksquare$

\bigskip\noindent\textbf{Theorem 4.5} \emph{Suppose} $p^{2}r<1$
\emph{and} $pr>1$\emph{. For every} $\alpha>0$ \emph{there exist}
$c_{8}=c_{8}(\alpha)>0$ \emph{and, with} $P_{p}$\emph{-probability
one, a finite} $\bar{k}_{8}(\eta)>0$ \emph{such that}

\begin{equation}
\sum_{i=1}^{m_{k}(\eta)}z_{i}^{2}(\eta)\leq c_{8}(1+\alpha)^{k}r^{k}\label{eq:3}\end{equation}
 \emph{for all} $k\geq\bar{k}_{8}(\eta)$.

\bigskip\noindent\textbf{Proof.} Let $\gamma=\frac{logr}{log(pr)}>1$
and $\gamma^{*}$ such that $\frac{1}{\gamma}+\frac{1}{\gamma^{*}}=1$
and take $\varepsilon_{1}$ such that $(1+\varepsilon_{1})^{\gamma^{*}}\leq(pr)^{1/3}$
and $(1+\varepsilon_{1})^{4}p^{2}r<1$. By Lemma 4.1, if $n\in\mathbb{N}$
and $V=V(n)\subseteq T^{(r)}$ is some set of vertices, then, since
$(1+\varepsilon_{1})^{\gamma^{*}}\leq(pr)^{1/3}$ we have

\begin{eqnarray*}
P_{p}(A_{V}(n)) & = & P_{p}(\exists v\in V(n):|C(v)\cap T_{n}^{(r)}|\geq(1+\varepsilon_{1})^{n-|v|}(pr)^{n-|v|})\\
 & \leq & \sum_{v\in V(n)}c_{5}e^{-c_{4}(1+\varepsilon_{1})^{\gamma^{*}(n-|v|)}}\end{eqnarray*}
\begin{equation}
\end{equation}
Recursively define $V_{j}$ and $d_{j}$ as follows:

\begin{eqnarray*}
V_{1}=V_{1}(n) & = & \left\{ v\in T^{(r)}:|v|\leq d_{1}n=n\frac{\log\left((1+\varepsilon_{1})^{4}p^{2}r\right)^{-1}}{\log r}\right\} ,\\
V_{j}=V_{j}(n) & = & \left\{ v\in T^{(r)},v\notin\bigcup_{j'=1}^{j-1}V_{j'}\right.\\
 &  & \left.:|v|\leq d_{j}n=n\frac{\log\left((1+\varepsilon_{1})^{4(1-d_{j-1})}p^{2(1-d_{j-1})}r^{1-2d_{j-1}}\right)^{-1}}{\log r}\right\} \end{eqnarray*}
\begin{equation}
\end{equation}
we then have

\begin{eqnarray*}
r^{d_{1}n} & = & \frac{1}{\left((1+\varepsilon_{1})^{4}p^{2}r\right)^{n}},\\
r^{d_{j}n} & = & \frac{1}{\left((1+\varepsilon_{1})^{4(1-d_{j-1})}p^{2(1-d_{j-1})}r^{1-2d_{j-1}}\right)^{n}},\end{eqnarray*}
\begin{equation}
\end{equation}

\begin{eqnarray*}
P_{p}(A_{V_{1}}(n)) & \leq & \left((1+\varepsilon_{1})^{4}p^{2}r\right)^{-n}c_{5}e^{-\left(\frac{1}{2}\right)^{\gamma^{*}}(1+\varepsilon_{1})^{\gamma^{*}n(1-d_{1})}},\\
P_{p}(A_{V_{j}}(n)) & \leq & \left((1+\varepsilon_{1})^{4(1-d_{j-1})}p^{2(1-d_{j-1})}r^{(1-2d_{j-1})}\right)^{-n}\\
 &  & \cdot c_{5}e^{-\left(\frac{1}{2}\right)^{\gamma^{*}}(1+\varepsilon_{1})^{\gamma^{*}n(1-d_{j})}}\end{eqnarray*}
\begin{equation}
\end{equation}
On $A_{V_{j}}(n)^{c}$ we have

\begin{eqnarray*}
\sum_{v\in V_{j}}|C(v)\cap T_{n}^{(r)}|^{2} & \leq & r^{d_{j}n}\left((1+\varepsilon_{1})pr\right)^{2n(1-d_{j-1})}\\
 & \leq & (1+\varepsilon_{1})^{-2n(1-d_{j-1})}r^{n}.\end{eqnarray*}
\begin{equation}
\end{equation}

Note that for $j=2,3,...$

\begin{equation}
d_{j}=(1-d_{j-1})\frac{\log(1+\varepsilon_{1})^{4}p^{2}r}{\log r}+d_{j-1}=(1-d_{j-1})d_{1}+d_{j-1}\label{eq:}\end{equation}
and that $d_{1}\in(0,1)$ since $(1+\varepsilon_{1})^{4}p^{2}r<1$,
so that $\lim_{j\rightarrow\infty}d_{j}=1$.

On the other hand, for the given $\alpha>0$ let $\rho_{1}$ be such
that $r^{\rho_{1}}<1+\alpha$; then, if for any cluster $C$ we let
$Base(C)=\min\{ k:C\cap T_{k}^{(r)}\neq\emptyset\}$, we have

\begin{eqnarray*}
\sum_{C:Base(C)\geq(1-\rho_{1})n}|C\cap T_{n}^{(r)}|^{2} & \leq & \sum_{C}|C\cap T_{n}^{(r)}|\max_{C:Base(C)\geq(1-\rho_{1})n}|C\cap T_{n}^{(r)}|\\
 & \leq & |T_{n}^{(r)}|r^{\rho_{1}n}\\
 & \leq & (1+\alpha)^{n}r^{n}.\end{eqnarray*}
\begin{equation}
\end{equation}

Next, take $J_{1}\in\mathbb{N}$ such that $d_{J_{1}}\geq(1-\rho_{1})$.
Then 

\begin{eqnarray*}
\sum_{n=1}^{\infty}\sum_{j=1}^{J_{1}}P_{p}(A_{V_{j}}(n)) & \leq & \sum_{j=1}^{J_{1}}\sum_{n=1}^{\infty}\left((1+\varepsilon_{1})^{4(1-d_{j-1})}p^{2(1-d_{j-1})}r^{(1-2d_{j-1})}\right)^{-n}\\
 &  & \cdot c_{5}e^{-c_{6}(1+\varepsilon_{1})^{\gamma^{*}n(1-d_{j})}}<+\infty\end{eqnarray*}
\begin{equation}
\end{equation}
since for each $j$ the series is of the form $A^{n}e^{-B^{n}}$,
with $A>1$ and $B>0$, thus convergent. This implies that, by Borel-Cantelli,
$A_{V_{1}}(n)\cup A_{V_{2}}(n)\cup...\cup A_{V_{J_{1},}}(n)$ occurs
only for a finite number of $n$'s with probability one. Thus, for
almost all $\eta$ there exists $\bar{k}_{8}(\eta)$ such that for
all $k>\bar{k}_{8}(\eta)$, $\bigcap_{j=1}^{J_{1}}A_{V_{j}}(k)^{c}$
occurs and this implies 

\begin{eqnarray*}
\sum_{i=1}^{m_{k}(\eta)}z_{i}^{2} & = & \sum_{C}|C\cap T_{k}^{(r)}|^{2}\\
 & \leq & \sum_{C:Base(C)\geq(1-\rho_{1})k}|C\cap T_{k}^{(r)}|^{2}+\sum_{j=1}^{J_{1}}\sum_{C:Base(C)\in V_{j}}|C(v)\cap T_{k}^{(r)}|^{2}\\
 & \leq & (1+\alpha)^{k}r^{k}+(1+\varepsilon_{1})^{-2k(1-d_{J_{1}})}r^{k}J_{1}\\
 & \leq & c_{8}(1+\alpha)^{k}r^{k}\end{eqnarray*}
\begin{equation}
\end{equation}
for a suitable $c_{8}=c_{8}(J_{1})$. $\blacksquare$

\bigskip\noindent\textbf{Theorem 4.6} \emph{If} $p^{2}r<1$ \emph{and}
$pr>1$\emph{, then there exist} $\bar{\alpha}'>0$, $c_{9}>0$ \emph{and,
with} $P_{p}$\emph{-probability one, a finite} $\bar{k}_{9}(\eta)>0$
\emph{such that for every $\alpha'<\bar{\alpha}'$}

\begin{equation}
\sum_{i=1}^{m_{k}(\eta)}z_{i}^{3}\leq c_{9}(1-\alpha')^{k}(pr^{2})^{k}\label{eq:4}\end{equation}
\emph{for all} $k\geq\bar{k}_{39}(\eta)$.

\bigskip\noindent\textbf{Proof.} We proceed as in the proof of Theorem
4.2 by taking $\varepsilon_{1}$, $V_{j}$, $A_{V_{j}}(n)$. On $A_{V_{j}}(n)^{c}$
we now have \begin{eqnarray*}
\sum_{v\in V_{j}}|C(v)\cap T_{n}^{(r)}|^{3} & \leq & r^{d_{j}n}((1+\varepsilon_{1})pr)^{3n(1-d_{j-1})}\\
 & \leq & (1+\varepsilon_{1})^{-n(1-d_{j.-1})}p^{n(1-d_{j-1})}r^{(2-d_{j-1})n}\\
 & \leq & \frac{1}{((1+\varepsilon_{1})^{1-d_{j-1}}(pr)^{d_{j-1}})^{n}}(pr^{2})^{n}\end{eqnarray*}
\begin{equation}
\end{equation}
with $d_{j}$'s defined as above.

Now, take $\rho_{2}>0$ such that $\rho_{2}<\frac{log(pr)}{4logr}$.
Then\begin{eqnarray*}
\sum_{C:Base(C)\geq(1-\rho_{2})n}|C\cap T_{n}^{(r)}|^{3} & \leq & \sum_{C}|C\cap T_{n}^{(r)}|\max_{C:Base(C)\geq(1-\rho_{2})n}|C\cap T_{n}^{(r)}|^{2}\\
 & \leq & r^{n}r^{2\rho_{2}n}\\
 & \leq & r^{n}(pr)^{n}(pr)^{-n/2}\\
 & \leq & (1-\alpha')^{n}(pr^{2})^{n}\end{eqnarray*}
\begin{equation}
\end{equation}
provided that $1-\alpha'\geq\frac{1}{\sqrt{pr}}$.

Next, take $J_{2}\in\mathbb{N}$ such that $d_{J_{2}}\geq1-\rho_{2}$
and note that the Borel-Cantelli Lemma applies as above. Take $\alpha'$
also satisfying $1-\alpha'\geq(1+\varepsilon_{1})^{-(1-d_{J_{2}})}$.
Then, for $k\geq\bar{k}_{9}(\eta)$, 

\begin{eqnarray*}
\sum_{i=1}^{m_{k}(\eta)}z_{i}^{3} & = & \sum_{C}|C\cap T_{k}^{(r)}|^{3}\\
 & \leq & \sum_{C:Base(C)\geq(1-\rho_{2})k}|C\cap T_{k}^{(r)}|^{3}+\sum_{j=1}^{J_{2}}\sum_{C:Base(C)\in V_{j}}|C(v)\cap T_{k}^{(r)}|^{3}\\
 & \leq & (1-\alpha')^{k}(pr^{2})^{k}+\frac{1}{(1+\varepsilon_{1})^{(1-d_{J_{2}})k}}(pr^{2})^{k}\\
 & \leq & c_{9}(1-\alpha')^{k}(pr^{2})^{k}.\end{eqnarray*}
\begin{equation}
\end{equation}
$\blacksquare$

The next result gives the inequality for critical points $p_{c}(k)$.

\bigskip         For the lower bound we need 

\bigskip\noindent\textbf{Lemma 4.8} \emph{If $Z_{i}$'s, $i=1,...,m$,
are independent random variables each taking value in some $\{-l,l\}$,
$l\in\mathbb{N}$ such that $Z_{i}\in\{-1,1\}$ for all $i=1,...,I$
then for every $\alpha>0$ and $m\geq I>0$ we have }

\begin{equation}
P(\sum_{i=1}^{m}Z_{i}\in[-\alpha,\alpha])\leq P(\sum_{i=1}^{I}Z_{i}\in[-\alpha,\alpha])\end{equation}

\bigskip\noindent\textbf{Proof.} Since $p_{k}=P(\sum_{i=1}^{I}Z_{i}=k)=\left(\begin{array}{c}
I\\
(I+k)/2\end{array}\right)2^{-I}$, $p_{k}$ increases up to $I/2$ and decreases afterwards; then,
letting $S_{k}=\sum_{i=1}^{k}Z_{i}$, we have 

\begin{eqnarray}
P(S_{m}\in[-\alpha,\alpha]) & = & P(S_{I}\in[-\alpha,\alpha],S_{m}\in[-\alpha,\alpha])\nonumber \\
 &  & +P(S_{I}\notin[-\alpha,\alpha],S_{m}\in[-\alpha,\alpha])\nonumber \\
 & = & P(S_{I}\in[-\alpha,\alpha],S_{m}\in[-\alpha,\alpha])\nonumber \\
 &  & +\sum_{t\notin[-\alpha,\alpha]}\sum_{l\in[-\alpha-t,\alpha-t]}P(S_{I}=t,S_{m-I}=l)\nonumber \\
 & \leq & P(S_{I}\in[-\alpha,\alpha],S_{m}\in[-\alpha,\alpha])\\
 &  & +\sum_{t\notin[-\alpha,\alpha]}\sum_{l\in[-\alpha-t,\alpha-t]}P(S_{I}=t+l,S_{m-I}=-l)\nonumber \\
 & = & P(S_{I}\in[-\alpha,\alpha],S_{m}\in[-\alpha,\alpha])\nonumber \\
 &  & +P(S_{I}\in[-\alpha,\alpha],S_{m}\notin[-\alpha,\alpha])\nonumber \\
 & = & P(S_{I}\in[-\alpha,\alpha])\nonumber \end{eqnarray}
\[
\quad\quad\quad\quad\quad\quad\quad\quad\quad\quad\quad\quad\quad\quad\quad\quad\quad\quad\quad\quad\quad\quad\quad\quad\blacksquare\]

For the upper bound we need an estimate for the error rate $\varepsilon(k)$
at distance $k$, i.e. the value defined by\begin{equation}
1-\varepsilon(k)=P_{p}\left(\sum_{v\in T_{k}^{(r)}}\sigma_{v}>0\left|\sigma_{0}=1\right.\right)+\frac{1}{2}P_{p}\left(\sum_{v\in T_{k}^{(r)}}\sigma_{v}=0\left|\sigma_{0}=1\right.\right)\label{eq:}\end{equation}

\bigskip\noindent\textbf{Lemma 4.7} \emph{If} $p^{2}r<1$ \emph{and}
$pr>1$ \emph{then there exists} $c_{10}>0$ \emph{such that for every}
$\alpha>0$ \emph{with probability one there exists} $\bar{k}_{11}$
\emph{finite such} \emph{that for all} $k>\bar{k}_{10}$ 

\begin{equation}
1-\varepsilon(k)\geq\frac{1}{2}+\frac{1}{2}c_{10}\frac{(p\sqrt{r})^{k}}{(1+\alpha)^{k/2}}\label{eq:}\end{equation}

\bigskip\noindent\textbf{Proof.} We have\begin{gather*}
P_{p}\left(\left.\sum_{v\in T_{k}^{(r)}}\sigma_{v}>0\right|\sigma_{0}=1\right)+\frac{1}{2}P_{p}\left(\left.\sum_{v\in T_{k}^{(r)}}\sigma_{v}=0\right|\sigma_{0}=1\right)\end{gather*}

\begin{gather}
=\left[P_{p}\left(\left.\sum_{v\in T_{k}^{(r)}}\sigma_{v}>0\right|\sigma_{0}=1,|R_{k}|<\frac{(pr)^{k}}{2}\right)\right.\label{eq:11}\\
\left.+\frac{1}{2}P_{p}\left(\left.\sum_{v\in T_{k}^{(r)}}\sigma_{v}=0\right|\sigma_{0}=1,|R_{k}|<\frac{(pr)^{k}}{2}\right)\right]\cdot P_{p}\left(\left.|R_{k}|<\frac{(pr)^{k}}{2}\right|\sigma_{0}=1\right)\nonumber \\
+\left[P_{p}\left(\left.\sum_{v\in T_{k}^{(r)}}\sigma_{v}>0\right|\sigma_{0}=1,|R_{k}|\geq\frac{(pr)^{k}}{2}\right)\right.\nonumber \\
\left.+\frac{1}{2}P_{p}\left(\left.\sum_{v\in T_{k}^{(r)}}\sigma_{v}=0\right|\sigma_{0}=1,|R_{k}|\geq\frac{(pr)^{k}}{2}\right)\right]\cdot P_{p}\left(\left.|R_{k}|\geq\frac{(pr)^{k}}{2}\right|\sigma_{0}=1\right)\nonumber \end{gather}

Notice that for each $\eta\in\{-1,1\}^{\mathcal{E}(T^{(r)})}$, $\sum_{v\in T_{k}^{(r)}}\sigma_{v}=\sum_{i=1}^{m_{k}(\eta)}Z_{i}+|R_{k}|$,
with $Z_{i}$ symmetric random variables. Therefore,\begin{gather}
P_{p}(\sum_{v\in T_{k}^{(r)}}\sigma_{v}>0|\sigma_{0}=1,|R_{k}|<\frac{(pr)^{k}}{2})+\frac{1}{2}P_{p}(\sum_{v\in T_{k}^{(r)}}\sigma_{v}=0|\sigma_{0}=1,|R_{k}|<\frac{(pr)^{k}}{2})\nonumber \\
\geq\sum_{\begin{array}{c}
\eta\in\{-1,1\}^{\mathcal{E}(T^{(r)})}\\
|R_{k}|\leq\frac{(pr)^{k}}{2}\end{array}}\left[P_{p}\left(\sum_{i=1}^{m_{k}(\eta)}Z_{i}>0\left|\eta\right.\right)+\frac{1}{2}P_{p}\left(\sum_{i=1}^{m_{k}(\eta)}Z_{i}=0\left|\eta\right.\right)\right]P_{p}(\eta)\geq\frac{1}{2}.\end{gather}

For the second part of (\ref{eq:11}) we use that

\begin{gather}
P_{p}(\sum_{v\in T_{k}^{(r)}}\sigma_{v}>0|\sigma_{0}=1,|R_{k}|\geq\frac{(pr)^{k}}{2})+\frac{1}{2}P_{p}(\sum_{v\in T_{k}^{(r)}}\sigma_{v}=0|\sigma_{0}=1,|R_{k}|\geq\frac{(pr)^{k}}{2})\nonumber \\
\geq\sum_{\eta:|R_{k}(\eta)|\geq\frac{(pr)^{k}}{2}}\left[P_{p}\left(\sum_{i=1}^{m_{k}(\eta)}Z_{i}>0\left|\eta\right.\right)+\frac{1}{2}P_{P}\left(\sum_{i=1}^{m_{k}(\eta)}Z_{i}=0\left|\eta\right.\right)\right.\nonumber \\
\left.+P_{p}\left(\sum_{i=1}^{m_{k}(\eta)}Z_{i}\in\left(-\frac{(pr)^{k}}{2},0\right]\left|\eta\right.\right)\right]\frac{P_{p}(\eta)}{P_{p}(|R_{k}|\geq\frac{(pr)^{k}}{2})}\nonumber \\
\geq\frac{1}{2}+\frac{1}{2}\sum_{\eta:|R_{k}(\eta)|\geq\frac{(pr)^{k}}{2}}P_{p}\left(\left|\sum_{i=1}^{m_{k}(\eta)}Z_{i}\right|<\frac{(pr)^{k}}{2}\left|\eta\right.\right)\frac{P_{p}(\eta)}{P_{p}(|R_{k}|\geq\frac{(pr)^{k}}{2})}\end{gather}
Then\begin{gather}
P_{p}(\sum_{v\in T_{k}^{(r)}}\sigma_{v}>0|\sigma_{0}=1)+\frac{1}{2}P_{p}(\sum_{v\in T_{k}^{(r)}}\sigma_{v}=0|\sigma_{0}=1)\nonumber \\
\geq\frac{1}{2}+\frac{1}{2}\sum_{\eta:|R_{k}(\eta)|\geq\frac{(pr)^{k}}{2}}P_{p}\left(\left|\sum_{i=1}^{m_{k}(\eta)}Z_{i}\right|\leq\frac{(pr)^{k}}{2}\left|\eta\right.\right)P_{p}(\eta)\label{eq:7}\end{gather}

Since the random variable $W$ defined in Lemma 4.1 is absolutely
continuous and $E(W)=1$ (see \cite{key-6}), then $P(W\geq\frac{1}{2})>0$.
Moreover, $\frac{|R_{k}|}{(pr)^{k}}$ converges in distribution to
$W$, so there exists a non random $\bar{k}_{1}$ such that for all
$k\geq\bar{k}_{1}$

\begin{equation}
P_{p}\left(\frac{|R_{k}|}{(pr)^{k}}\geq\frac{1}{2}\right)\geq\frac{1}{2}P(W\geq\frac{1}{2})>0.\end{equation}

We then want to estimate $P_{p}\left(\left|\sum_{i=1}^{m_{k}(\eta)}Z_{i}\right|\leq\frac{(pr)^{k}}{2}\left|\eta\right.\right)$
via the Gaussian approximation using the Berry-Essen estimates of
the error. To this extent, we will use the results in Theorems 4.2,
4.3 and 4.4 with $\alpha$ of Theorem 4.3 such that $(1+\alpha)^{-1/2}>1-\alpha'$,
with $\alpha'<\bar{\alpha}'$ and $\bar{\alpha}'$ determined as in
Theorem 4.4. Such results hold with $P_{p}$-probability one for almost
all $\eta$'s, and thus it is possible to find a non random $\bar{k}_{2}$
such that $P_{p}(\eta:\bar{k}_{2}\geq max(\bar{k}_{7}(\eta),\bar{k}_{8}(\eta),\bar{k}_{9}(\eta))>1-\frac{1}{4}P(W\geq\frac{1}{2})$.
Let $\bar{k}_{3}$ such that $\left(\frac{p^{2}r}{1+\alpha}\right)^{k}\frac{1}{4c_{8}(\alpha)}<-\log\frac{1}{2}$
and $\frac{1}{\sqrt{c_{8}}(1+\alpha)^{k}}\geq2\frac{c_{9}}{c_{7}^{3/2}}(1-\alpha')^{k}$,
for $k>\bar{k}_{3}$. 

If we define the non random constant

\begin{equation}
\bar{k}_{10}=max(\bar{k}_{1},\bar{k}_{2},\bar{k}_{3})\end{equation}

and 

\begin{gather}
M_{k}=\left\{ \eta\in\{-1,1\}^{\mathcal{E}(T^{(r)})}\left|\frac{|R_{k}(\eta)|}{(pr)^{k}}\geq\frac{1}{2},c_{7}r^{k}\leq\sum_{i=1}^{m_{k}(\eta)}Z_{i}^{2}(\eta)\leq c_{8}(1+\alpha)^{k}r^{k},\right.\right.\nonumber \\
\left.\sum_{i=1}^{m_{k}(\eta)}|Z_{i}^{3}(\eta)|\leq c_{9}(1-\alpha')^{k}(pr^{2})^{k}\right\} \end{gather}

then, for $k\geq\bar{k}_{11}$

\begin{equation}
P_{p}(M_{k})\geq\frac{1}{4}P(W\geq\frac{1}{2})>0.\label{eq:8}\end{equation}

From (\ref{eq:7}) we then get

\begin{gather}
P_{p}(\sum_{v\in T_{k}^{(r)}}\sigma_{v}>0|\sigma_{0}=1)+\frac{1}{2}P_{p}(\sum_{v\in T_{k}^{(r)}}\sigma_{v}=0|\sigma_{0}=1)\nonumber \\
\geq\frac{1}{2}+\frac{1}{2}\sum_{\eta\in M_{k}}P_{p}\left(\left|\sum_{i=1}^{m_{k}(\eta)}Z_{i}\right|\leq\frac{(pr)^{k}}{2}\left|\eta\right.\right)P_{p}(\eta),\label{eq:9}\end{gather}
which we now estimate using the Gaussian approximation. Given $\eta$,
the $Z'_{i}$'s are independent random variables, so we can substitute
them with the equally distributed $Z'_{i}$'s. The Berry-Essen Theorem
gives

\begin{eqnarray*}
P\left(\sum_{i=1}^{m_{k}(\eta)}Z'_{i}\in\left[-\frac{(pr)^{k}}{2},\frac{(pr)^{k}}{2}\right]\right) & = & P\left(\sum_{i=1}^{m_{k}(\eta)}\frac{Z_{i}}{\sqrt{V_{k}}}\in\left[\frac{-\frac{(pr)^{k}}{2}}{\sqrt{V_{k}}},\frac{\frac{(pr)^{k}}{2}}{\sqrt{V_{k}}}\right]\right)\\
 & = & \int_{-\frac{\frac{(pr)^{k}}{2}}{\sqrt{V_{k}}}}^{\frac{\frac{(pr)^{k}}{2}}{\sqrt{V_{k}}}}\frac{1}{\sqrt{2\pi}}e^{-x^{2}/2}dx+E_{k}\end{eqnarray*}
\begin{equation}
\end{equation}
with $|E_{k}|\leq\frac{s_{k}}{V_{k}^{3/2}}$, where $V_{k}=\sum_{i=1}^{m_{k}(\eta)}Var(Z_{i})=\sum_{i=1}^{m_{k}(\eta)}z_{i}^{2}$
and $s_{k}=\sum_{i=1}^{m_{k}(\eta)}E(|z_{i}|^{3})=\sum_{i=1}^{m_{k}}z_{i}^{3}$.

If $\eta\in M_{k}$, $V_{k}\leq c_{8}(1+\alpha)^{k}r^{k}$ and

\begin{equation}
E_{k}\leq\frac{c_{9}(1-\alpha')^{k}(pr^{2})^{k}}{(c_{1}r^{k})^{3/2}}=\frac{c_{9}}{c_{7}^{3/2}}(1-\alpha')^{k}p^{k}r^{k/2}\end{equation}
so that 

\begin{eqnarray*}
P\left(\left|\sum_{i=1}^{m_{k}(\eta)}Z_{i}\right|\leq\frac{(pr)^{k}}{2}\right) & \geq & \int_{-\frac{\frac{1}{2}(pr)^{k}}{\sqrt{c_{8}(1+\alpha)^{k}r^{k}}}}^{\frac{\frac{1}{2}(pr)^{k}}{\sqrt{c_{8}(1+\alpha)^{k}r^{k}}}}\frac{1}{\sqrt{2\pi}}e^{-x^{2}/2}dx-\frac{c_{9}}{c_{7}^{3/2}}(1-\alpha')^{k}p^{k}r^{k/2}\\
 & \geq & \frac{p^{k}r^{k/2}}{\sqrt{c_{8}}(1+\alpha)^{k/2}}e^{-\frac{p^{2k}r^{k}}{4c_{7}(1+\alpha)^{k}}}-\frac{c_{9}}{c_{7}^{3/2}}(1-\alpha')^{k}p^{k}r^{k/2}\\
 & \geq & \frac{1}{2}\frac{p^{k}r^{k/2}}{\sqrt{c_{8}}(1+\alpha)^{k/2}}\end{eqnarray*}
\begin{equation}
\end{equation}
for $k\geq\bar{k}_{10}\geq\bar{k}_{3}$. Together with (\ref{eq:8}),
(\ref{eq:9}) this implies 

\begin{equation}
P(\sum_{v\in T_{k}^{(r)}}\sigma_{v}>0|\sigma_{0}=1)+\frac{1}{2}P(\sum_{v\in T_{k}^{(r)}}\sigma_{v}=0|\sigma_{0}=1)\geq\frac{1}{2}+\frac{1}{2}c_{10}\frac{p^{k}r^{k/2}}{(1+\alpha)^{k/2}}\end{equation}
with $c_{10}=\frac{1}{\sqrt{c_{8}}}P_{p}(M_{k})>0$, for all $k\geq\bar{k}_{10}$.
$\blacksquare$

\bigskip\noindent\textbf{Proof of Theorem 4.1} From Lemma 4.7, the
probability of error free transmission $p(k)=1-2\varepsilon(k)$ satisfies 

\begin{equation}
p(k)\geq c_{11}\frac{(p\sqrt{r})^{k}}{(1+\alpha)^{k/2}}\label{eq:new4}\end{equation}
for the binary transmission problem on $T^{(r^{k})}$for $k$ large
enough. Therefore, there is reconstruction if

\begin{equation}
1<p(k)r^{k/2}=c_{11}\left(\frac{pr}{\sqrt{1+\alpha}}\right)^{k},\end{equation}
which is to say

\begin{equation}
p_{c}(k)\leq\frac{1+\alpha}{c_{11}^{1/k}r}\label{eq:new3}\end{equation}
for $k$ large enough. Let $\alpha_{k}$ be the smallest $\alpha$
s.t. (\ref{eq:new3}) holds. Then $\lim_{k\rightarrow\infty}\alpha_{k}=0$
as required to prove the upper bound of Theorem 4.1.

\medskip

Similarly to (\ref{eq:7}) we estimate, for $\beta>0$, 

\begin{eqnarray*}
1-\varepsilon(k) & = & \left[P_{p}\left(\left.\sum_{v\in T_{k}^{(r)}}\sigma_{v}>0\right|\sigma_{0}=1,|R_{k}|>(1+\beta)^{k}(pr)^{k}\right)\right.\\
 &  & \left.+\frac{1}{2}P_{p}\left(\left.\sum_{v\in T_{k}^{(r)}}\sigma_{v}=0\right|\sigma_{0}=1,|R_{k}|>(1+\beta)^{k}(pr)^{k}\right)\right]\\
 &  & \cdot P_{p}(|R_{k}|>(1+\beta)^{k}(pr)^{k}|\sigma_{0}=1)\\
 &  & +\left[P_{p}\left(\left.\sum_{v\in T_{k}^{(r)}}\sigma_{v}>0\right|\sigma_{0}=1,|R_{k}|\leq(1+\beta)^{k}(pr)^{k},\left|\sum_{i=1}^{m_{k}(\eta)}Z_{i}\right|\leq(1+\beta)^{k}(pr)^{k}\right)\right.\\
 &  & \left.+\frac{1}{2}P_{p}\left(\left.\sum_{v\in T_{k}^{(r)}}\sigma_{v}=0\right|\sigma_{0}=1,|R_{k}|\leq(1+\beta)^{k}(pr)^{k},\left|\sum_{i=1}^{m_{k}(\eta)}Z_{i}\right|\leq(1+\beta)^{k}(pr)^{k}\right)\right]\\
 &  & \cdot P_{p}(|R_{k}|\leq(1+\beta)^{k}(pr)^{k},\left|\sum_{i=1}^{m_{k}(\eta)}Z_{i}\right|\leq(1+\beta)^{k}(pr)^{k}|\sigma_{0}=1)\\
 &  & +\left[P_{p}\left(\left.\sum_{v\in T_{k}^{(r)}}\sigma_{v}>0\right|\sigma_{0}=1,|R_{k}|\leq(1+\beta)^{k}(pr)^{k},\left|\sum_{i=1}^{m_{k}(\eta)}Z_{i}\right|>(1+\beta)^{k}(pr)^{k}\right)\right.\\
 &  & \left.+\frac{1}{2}P_{p}\left(\left.\sum_{v\in T_{k}^{(r)}}\sigma_{v}=0\right|\sigma_{0}=1,|R_{k}|\leq(1+\beta)^{k}(pr)^{k},\left|\sum_{i=1}^{m_{k}(\eta)}Z_{i}\right|>(1+\beta)^{k}(pr)^{k}\right)\right]\\
 &  & \cdot P_{p}(|R_{k}|\leq(1+\beta)^{k}(pr)^{k},\left|\sum_{i=1}^{m_{k}(\eta)}Z_{i}\right|>(1+\beta)^{k}(pr)^{k}|\sigma_{0}=1).\end{eqnarray*}

\begin{equation}
\label{eq:newnew}\end{equation}

From Lemma 4.3 \[
P_{p}(|R_{k}|>(1+\beta)^{k}(pr)^{k}|\sigma_{0}=1)\leq c_{5}e^{c_{6}(1+\beta)^{\gamma^{*}k}}.\]

In the third term, the expression between square brackets is exactly
$\frac{1}{2}$ since $Z_{i}$'s are independent and symmetric.

Next we consider the second term. Assume first $pr\geq1$. Let $I$
be the set of vertices in $T_{k}^{(r)}$ which are isolated FK clusters.
Then, by large deviations for i.i.d. random variables, $P_{p}(|I|<\frac{1}{2}\left(\frac{1-p}{r}\right)r^{k}\leq e^{-cr^{k}}$.
Moreover, from Lemma 4.8, the expression between square brackets in
the second term of (\ref{eq:newnew}) is bounded above by $P_{p}\left(\left|\sum_{i=1}^{I}Z'_{i}\right|\leq((1+\beta)pr)^{k}\right)$,
with $Z'_{i}$ i.i.d. symmetric random variables with values in $\{-1,1\}$.
In turn, if $|I|\geq(1-p-\varepsilon)r^{k}\geq\frac{r^{k}}{2}$ the
normal approximation implies that for some $c_{12}>0$, $c_{13}>0$,

\begin{eqnarray*}
P_{p}\left(\left|\sum_{i=1}^{I}Z'_{i}\right|\leq((1+\beta)pr)^{k}\right) & \leq & c_{12}((1+\beta)p\sqrt{r})^{k}+c_{13}\frac{1}{\sqrt{I}}\\
 & \leq & c_{1}p^{k}r^{k/2}(1+\beta)^{k}\end{eqnarray*}

\begin{equation}
\end{equation}
for a suitable $c_{1}$ large enough, where the last term comes from
the Berry-Essen error estimate for the random variables $Z'_{i}$,
with $|I|\geq\frac{r^{k}}{2}$ and $\frac{1}{r^{k/2}}\leq p^{k}r^{k/2}$.

Collecting the above estimates we have

\begin{eqnarray*}
1-\varepsilon(k) & \leq & e^{-cr^{k}}+c_{12}((1+\beta)p\sqrt{r})^{k}+c_{13}p^{k}r^{k/2}(1-\alpha')^{k}+\frac{1}{2}\\
 & \leq & \frac{1}{2}+\frac{1}{2}c_{1}p^{k}r^{k/2}(1+\beta)^{k}.\end{eqnarray*}

\begin{equation}
\end{equation}

Therefore,

\begin{eqnarray}
p(k) & \leq & c_{1}p^{k}r^{k/2}(1+\beta)^{k}\label{eq:newnew4}\end{eqnarray}
and the condition for non-reconstruction on the rescaled tree $T^{(r^{k})}$
becomes\begin{equation}
c_{1}p^{k}r^{k}(1+\beta)^{k}<1.\end{equation}

This implies

\begin{equation}
p_{c}(k)\geq\frac{1}{(1+\beta)c_{1}^{1/k}r}\geq\frac{1}{c_{1}^{1/k}r}.\end{equation}

If, on the other hand, $pr<1$, then for small enough $\beta$, $(1+\beta)pr<1$
and the second term in square brackets of (\ref{eq:newnew}) reduces
to $\frac{1}{2}P_{p}\left(\left.\sum_{v\in T_{k}^{(r)}}\sigma_{v}=0\right|\sigma_{0}=1\right)$,
but clearly in this case the symmetry is not broken and no reconstruction
can take place. $\blacksquare$

From Theorem 4.1 it is obvious that the critical points $p_{c}(k)$
converge to the Ising model critical point.

\section{Minority removal}

The self-correction mechanism discussed above is not suitable for
biological transmission, in which offsprings, once generated, cannot
be changed. However, there is a similar mechanism, which consists
of self-correcting a generation by removing the elements not belonging
to the majority, which could be implemented in a biological setting.
If $r\geq4$ and such minority removal is carried out every step in
blocks of size $M$, then in the renormalized tree each (macroscopic)
vertex has a random number of children larger then or equal to $2$,
while the error rate is estimated as in (\ref{eq:2.4}) but on a random
number of vertices, between $\frac{M}{2}$ and $M$; by taking inequalities
as done below, one can see that (\ref{eq:2.5}) still holds with minor
changes and thus reconstruction is also possible at every $\varepsilon<\frac{1}{2}$
with a sufficiently large $M$. It is also the case that if a within-descent
minority removal is carried out every $k$ generations, only minor
changes in the constants are needed in Theorem 4.1 and the limit of
the critical points is still the Ising critical point as in Corollary
4.2.

This highlights a possibly real but rather particular phenomenon.
It looks like a bit of information in the parent biological unit is
better transmitted, i.e. it is more resistant to random transmission
errors, if enhanced by regularly destroying descendants not belonging
to the local majority. From the biological point of view this is also
likely to improve the functionality of local segments (cells or individuals,
for instance). However, the minority removal sometimes preserves the
wrong information, thus creating blocks of mutated descendants, a
phenomenon similar to tumor formation. In this respect, our findings
seem to suggest that tumor generation might be intrinsically connected
to improvement in character transmission. Of course, any such claim
must be warranted by the study of many bits transmission.

Back to our single bit model, the minority-removal carried out every
step by blocks of size $M$ corresponds to first generating a random
tree $T_{M}'$ by means of a transformation $\Phi_{M}'$ analogous
to $\Phi_{M}$ and then identifying each block (of random size between
$\frac{M}{2}$ and $M$) by means of a transformation $\Psi_{M}'$,
analogous to $\Psi_{M}$. Let $\bar{P}_{\varepsilon,M}'=\Psi_{M}'(\Phi_{M}'(P_{\varepsilon}))$
be the distribution on the resulting random tree $T_{M}'$.

Similarly, the within-descent minority removal carried out every $k$-steps
corresponds to generating a random tree $T_{k}'$ by means of a transformation
$\Phi_{k}'$, analogous to $\Phi_{k}$, and then identifying each
block (of random size between $\frac{r^{k}}{2}$ and $r^{k}$) by
means of a transformation $\Psi_{k}'$, analogous to $\Psi_{k}$.
Let $P_{\varepsilon}'^{(k)}=\Psi_{k}'(\Phi_{k}'(P_{\varepsilon}))$
be the distribution on the resulting random tree $T_{k}'$.

Note that $T'_{M}$ and $T'_{k}$ are Galton-Watson trees, since they
are random trees with an i.i.d. number of offsprings in each vertex.
In generating $T'_{M}$ at least $M/2$ vertices are preserved in
each block of size $M$; these have at least $rM/2\geq2M$ descendants
which can be divided into at least $2$ blocks of size $M$ (and possibly
one remaining smaller block). Thus the number of descendants is at
least $2$. In generating $T'_{r}$ on the other hand, at least $r^{k}/2$
vertices are preserved in each block of size $r^{k}$ and each such
vertex gives rise to one descendant block, so each block (which is
a renormalized vertex) has at least $r^{k}/2$ (and at most $r^{k}$)
descendants.

The branching numbers, which on the Galton-Watson trees equal the
mean offspring number (see \cite{key-10}), satisfy then $br(T'_{M})\geq2$
and $r^{k}/2\leq br(T'_{r})\leq r^{k}$.

We begin with a Lemma stating that if on a subtree $T'\subset T$
maximum likelihood reconstruction takes place, then it does also on
$T$.

\bigskip\noindent  \textbf{Lemma 5.1} \emph{Given trees $T'\subseteq T$,
if maximum likelihood reconstruction takes place on $T'$ then it
does also on $T$, i.e. if $\liminf_{n}\Delta_{n}(P_{T'})>0$ then
$\liminf_{n}\Delta_{n}(P_{T})>0$.}

\bigskip\noindent  \textbf{Proof.} Let $A_{n}=\{\sigma_{n}\in T_{n}:P(\sigma_{n}|\sigma_{0}=+1)>P(\sigma_{n}|\sigma_{0}=-1)\}$,
let $A'_{n}$ be the same with $T_{n}$ replaced by $T'_{n}$ and
let $B'=\{\sigma'_{n}\in T'_{n}:P(\sigma'_{n}|\sigma_{0}=+1)=P(\sigma'_{n}|\sigma_{0}=-1)\}$.
We know $P(A'_{n}|\sigma_{0}=+1)-P(A'_{n}|\sigma_{0}=-1)\geq\delta>0$
for some $\delta$ for large $n$, and we want to show the same for
$A_{n}$. However, denoting by $P^{\pm}(\cdot)=P(\cdot|\pm1)$ we
have $P^{\pm}(A_{n}\cap(A'_{n})^{c})=P^{\mp}(A_{n}^{c}\cap A'_{n})$
by symmetry, and for any event $C$, by definition of $A_{n}$,

\begin{eqnarray}
P^{+}(A_{n}\cap C) & \geq & P^{-}(A_{n}\cap C)\nonumber \\
P^{+}(A_{n}^{c}\cap C) & \leq & P^{-}(A_{n}^{c}\cap C).\end{eqnarray}
Then,

\begin{gather*}
P^{+}(A_{n})-P^{-}(A_{n})\\
=P^{+}(A_{n}\cap A'_{n})+P^{+}(A_{n}\cap(A'_{n})^{c})+P^{+}(A_{n}\cap B')\\
-P^{-}(A_{n}\cap A'_{n})-P^{-}(A_{n}\cap(A'_{n})^{c})-P^{-}(A_{n}\cap B')\\
=P^{+}(A_{n}\cap A'_{n})+P^{-}(A_{n}^{c}\cap A'_{n})+P^{+}(A_{n}\cap B')\\
-P^{-}(A_{n}\cap A'_{n})-P^{+}(A_{n}^{c}\cap A'_{n})-P^{-}(A_{n}\cap B')\\
\geq P^{+}(A_{n}\cap A'_{n})+P^{+}(A_{n}^{c}\cap A'_{n})\\
-P^{-}(A_{n}\cap A'_{n})-P^{-}(A_{n}^{c}\cap A'_{n})\\
=P^{+}(A'_{n})-P^{-}(A'_{n})\end{gather*}

\begin{equation}
\end{equation}
from which the result follows. $\blacksquare$

The results for minority removal can be summarized as follows. Notice
that in the proof we use maximum likelihood reconstruction to use
Lemma 5.1 and get a bound on the critical point; on the other hand,
it is shown in \cite{key-1} that for binary tree the critical points
for majority or maximum likelihood reconstruction coincide. 

\bigskip\noindent  \textbf{Theorem 5.2 }

\begin{description}
\item [i)]\emph{If $r\geq4$, in the minority removal carried out every
step with blocks of size} $M$\emph{, for every} $\varepsilon\in[0,1/2)$
$\exists\bar{\bar{M}}:\forall M>\bar{\bar{M}}$ 
\end{description}
\[
\liminf_{n}\Delta_{n}(\bar{P}_{\varepsilon,M}')>0.\]

\begin{description}
\item [ii)]\emph{In the within-descent minority removal carried out every}
$k$ \emph{steps if} $p_{c}'(k)$ \emph{is the critical point then
with} $c>0$ \emph{as in Theorem 4.1 we have}
\end{description}
\[
\frac{1}{2^{\frac{1}{2k}}r}\leq p_{c}'(k)\leq\frac{4^{\frac{1}{2k}}}{c^{\frac{1}{2k}}r}\]

$\quad$\emph{so that }

\[
\lim_{k\rightarrow\infty}p_{c}'(k)=\frac{1}{r}.\]

\bigskip\noindent  \textbf{Proof.} i) In generating $T'_{M}$ at
least $\frac{M}{2}$ vertices were preserved in each block of size
$M$; these vertices have $r\frac{M}{2}\geq4\frac{M}{2}=2M$ descendants
which can be divided into at least two blocks of size $M$ (and some
remaining others, possibly smaller). Thus, the number of descendants
in the renormalized tree is at least $2$. 

On the other hand, the error rate $\bar{\varepsilon}_{M}'$ satisfies
(\ref{eq:2.5}) with $M$ replaced by $\frac{M}{2}$. By Lemma 5.1,
maximum likelihood reconstruction on $T'_{M}$ follows from that on
$T^{(M/2)}$ which is ensured by

\begin{equation}
2(1-2\bar{\varepsilon}_{M}')^{2}\geq2(1-2e^{-c_{\varepsilon}M/2})^{2}>1\end{equation}
which is satisfied for large $M$.

ii) In generating $T'_{r}$ at least $\frac{r^{k}}{2}$ vertices are
preserved in each block of size $r^{k}$; each such vertex gives rise
to one descendant block, so the branching number of the renormalized
tree is at least $\frac{r^{k}}{2}$. 

Also, it is possible to show bounds on the renormalized error free
transmission $p'(k)$ similar to those used to prove Theorem 4.1.
By carefully going through that proof, one can see that if $pr\geq1$

\begin{equation}
p'(k)^{2}\leq2c_{1}^{2}p^{2k}r^{k}(1+\beta)^{k}\end{equation}
as in (\ref{eq:newnew4}) if $pr<1$ again $p(k)$ is exponentially
small in $k$ and thus there is no reconstruction; and, finally 

\begin{equation}
p'(k)\geq c_{2}\frac{p^{k}\sqrt{\frac{r^{k}}{2}}}{(1+\alpha)^{k/2}}.\label{eq:}\end{equation}
as in (\ref{eq:new4}).

Again by Lemma 5.1 this implies

\[
\frac{1}{r}\vee\frac{1}{\left(2c_{1}\right)^{\frac{1}{2k}}r}\leq p_{c}'(k)\leq\frac{4^{\frac{1}{2k}}}{c_{2}^{\frac{1}{2k}}r}\]
 and 

\[
\lim_{k\rightarrow\infty}p_{c}'(k)=\frac{1}{r}.\]
 \[
\quad\quad\quad\quad\quad\quad\quad\quad\quad\quad\quad\quad\quad\quad\quad\quad\quad\quad\quad\quad\quad\quad\quad\quad\blacksquare\]

\textbf{Acknowledgments:} We thank G. Giacomin for useful discussions
and comments.

\lyxaddress{Alberto Gandolfi, Dipartimento di Matematica {}``Ulisse Dini'',
Viale Morgagni 67a, 50134 FIRENZE, ITALY, gandolfi@math.unifi.it}

\lyxaddress{Roberto Guenzani, Dipartimento di Matematica {}``Federico Enriques'',
Via Cesare Saldini, 50 MILANO ITALY, guenzani@mat.unimi.it}
\end{document}